%% file: unitarycasesubmitted.tex
\documentclass[leqno]{report}
\usepackage{amsmath,amsthm,amsfonts,amssymb,setspace,graphicx,
float, floatflt,wrapfig,verbatim}
\newcommand{\frontmatter}{\clearpage \pagenumbering{roman}}
\newcommand{\mainmatter}{\clearpage \pagenumbering{arabic}}

\input{artpreamble.tex}

\begin{document}
\pagestyle{headings}
\frontmatter
\begin{abstract}
\begin{center}
Stability of the local gamma factor in the unitary case\\
July  2006
\end{center}
In \cite{rallissoudry05}, Rallis and Soudry prove the stability
under twists by highly ramified characters
of the local gamma factor arising from the doubling method,
in the case of a symplectic group or orthogonal group $G$ over
a local non-archimedean field $F$ of characteristic zero,
and a representation $\pi$ of $G$, which is not necessarily
generic.  This paper extends their arguments to show
the stability in the case when $G$
is a unitary group over a quadratic extension $E$ of $F$,
thereby completing the proof of the stability for classical groups.
This stability property is important in Cogdell, Piatetski-Shapiro,
and Shahidi's use of the converse theorem to prove
the existence of a weak lift from automorphic, cuspidal,
generic representations of $G(\mathbb{A})$ to automorphic
representations of $\GL{n}{\mathbb{A}}$ for appropriate $n$,
to which references are given in \cite{rallissoudry05}.
\end{abstract}
\begin{titlepage}
\begin{singlespace}
\vspace*{50pt}
\begin{center} \Huge \bfseries
Stability of the local gamma factor in the unitary case.
\end{center}
\vspace{2.25in}
\begin{center}
A Paper
\end{center}
\vspace{2.25in}
\begin{center}
by \\ Eliot Brenner
\end{center}
\begin{center}
Center for Advanced Studies in Mathematics\\
at Ben Gurion University
\end{center}
\begin{center}
July 2006
\end{center}
\end{singlespace}
\end{titlepage}

\mainmatter
\newcounter{romnine}\setcounter{romnine}{9}
\section{Introduction}\label{sec:intro}
Let $G$ be either a symplectic or orthogonal group over
a local non-archimedean field $F$ of characteristic zero,
or a unitary group over a quadratic extension $E$ of $F$.
We consider the local gamma factor, associated to an irreducible
admissible representation $\pi$ of $G$, by the doubling
method of Piatetski-Shapiro and Rallis (\cite{gelbartpsrallis85},
\cite{lapidrallis}).  Denote the local gamma factor by
$\gamma(\pi,\chi,s,\psi)$, where $\chi$ is a character
of $F^*$, and $\psi$ is a fixed non-trivial character
$F$.  In this paper we treat the last open case of
the following result (cf. Theorem 1 in \cite{rallissoudry05}).
\begin{thm} \label{thm:generalstability} The local gamma factor
$\gamma(\pi,\chi,s,\psi)$ is stable, for $\chi$ sufficiently
ramified.  This means that for two given irreducible
admissible representations $\pi_1,\pi_2$ of $G$,
there exists an integer $N>0$, such that
\startdisp
\gamma(\pi_1,\chi,s,\psi)=\gamma(\pi_2,\chi,s,\psi),
\finishdisp
for all characters $\chi$, with conductor having an
exponent larger than $N$.
\end{thm}
Rallis and Soudry, in \cite{rallissoudry05}, prove
Theorem \ref{thm:generalstability} in the symplectic and
orthogonal cases, and this paper completes the proof
of Theorem \ref{thm:generalstability} by extending their arguments
to the unitary case.

The stability property of the local gamma factor, under
highly ramified twists, is well known for
$\mathrm{GL}_n\times\mathrm{GL}_m$.  It was proved
by Jacquet and Shalika.
For generic representation of split
classical groups, it is known thanks to the works of Cogdell,
Piatetski-Shapiro and Shahidi.  The stability
property is a key ingredient in the proof,
by the converse theorem, of the existence of a weak
lift from autmorphic, cuspidal, generic representations
of $G(\mathbb{A})$ ($G$ a split classical group)
to automorphic representations of
$\GL{n}{\mathbb{A}}$ (appropriate $n$), where
$\mathbb{A}$ is the adele ring of a given number field.
See \cite{rallissoudry05} for precise references to the literature.
In Theorem \ref{thm:generalstability}, $\pi$
is any irreducible representation of $G$; even when $G$
is quasi-split, $\pi$ is not necessarily generic.
The proof of the stability in this paper follows
the argument of \cite{rallissoudry05} closely.
Therefore, experts in the subject
will want to turn to the parts that are new and specific
to the unitary case; these are Lemmas \ref{lem:innerintegralvanishing}
and \ref{lem:innerintegralvanishingLzero},
some of the details of the calculations in
\ref{lem:cayleychofvar}, the standard
material contained in \S\S\ref{subsec:quadraticextgeneral}
through \ref{subsec:fieldnormextensions}, and, finally,
Proposition \ref{prop:stayingawayfrom1}.  The reader
who is familiar with \cite{rallissoudry05}
is advised to turn to Proposition \ref{prop:stayingawayfrom1}
first, since this elementary, but apparently new, observation
is the heart of the matter concerning the extension
of the arguments of \cite{rallissoudry05} to the
unitary case.

Recall that in the local theory of the doubling method,
we consider the integrals
\starteqn\label{eqn:zetaintegradefn}
Z(v_1,\hat{v}_2,f_{\chi,s})=
\int_G\langle\pi(g)v_1,\hat{v}_2\rangle
f_{\chi,s}(i(g,1))\intd g.
\finisheqn
Here $v_1$ lies in $V_{\pi}$---a space for $\pi$,
and $\hat{v}_2$ lies in the smooth dual of $V_{\pi}$,
$\hat{V}_{\pi}$ (affording the contragredient
representation $\hat{\pi}$).  Thus,
$g\mapsto \langle\pi(g)v_1,\hat{v}_2\rangle$
is a matrix coefficient of $\pi$; $f_{\chi, s}$
is a holomorphic section in an induced representation
of the split ``doubled" group $H$---induced from the Siegel
parabolic subgroup $P$ of $H$, and a character,
which is of the form
$\chi(\det\cdot)|\det\cdot|^{s-1/2}$.  Finally, there is
an embedding $i: G\times G\rightarrow H$, such that
$P\cdot i(G\times G)=P\cdot i(G\times 1)$ is an
open and dense subset in $H$.  The integrals
\eqref{eqn:zetaintegradefn} converge absolutely in a right-half-plane
and continue meromorphically to the whole plane,
being rational functions in $q^{-s}$, where $q$ is
the number of elements in the residue field of $F$.
The functions $Z(v_1,\hat{v}_2,f_{\chi,s})$
satisfy a functional equation
\startdisp
\Gamma(\pi,\chi,s)Z(v_1,\hat{v}_2,f_{\chi,s})=
Z(v_1,\hat{v}_2,M(\chi,s)f_{\chi,s}),
\finishdisp
where $M(\chi,s)$ is the intertwining operator
associated to the element $w=i(1,-1)$.
The proportionality factor $\Gamma(\pi,\chi,s)$
is a rational function of $q^{-s}$
which
depends only on $\pi$ and $\chi$.  Note that
$\Gamma(\pi,\chi,s)$ is independent of $\psi$.

The local gamma factor $\gamma(\pi,\chi,s,\psi)$
is obtained from
$\omega_{\pi}(-1)\Gamma(\pi,\chi,s)$ (where $\omega_{\pi}$
is the central character of $\pi$) by multiplication
by a factor which depends on $\chi,\psi$ (and $G$)
and not on $\pi$.  See \cite{rallissoudry05},
pp. 292--3 for the details.  Therefore, Theorem
\ref{thm:generalstability} will follow from
\begin{thm}  Let $\pi$ be an irreducible admissible
representation of $G$.  Then $\omega_{\pi}(-1)
\Gamma(\pi,\chi,s)$ is stable, for sufficiently
ramified $\chi$.  More precisely, there is a positive
integer $N$, such that for all ramified characters $\chi$
of $F^*$, with conductor having exponent larger than
$N$, we have
\startdisp
\omega_{\pi}(-1)\Gamma(\pi,\chi,s)=M(\chi,s)f_{\chi,s}
(i(-1,1))
\finishdisp
for certain choice of $f_{\chi,s}$.
\end{thm}
A more precise form of this stability is given in Theorem
\ref{thm:main}.
\vspace*{0.3cm}

\noindent\textbf{Acknowledgements.}\hspace*{0.5mm}
The author thanks David Soudry for suggesting this problem
and thanks Soudry, Nadya Gurevich, and Omer Offen for helpful conversations.
He also thanks Mr.\hspace*{-.5mm} Tony Petrello for additional financial
support during the writing of this paper.

\section{Notation and Preliminaries}\label{sec:notandprelim}
As far as possible, we keep the notation consistent with
\cite{rallissoudry05}.
Let $E$ be any local non-archimedean field, of characteristic zero.
We denote by $\scrO_E$ its ring of integers, and by $\scrP_E$
the prime ideal of $\scrO_E$.  We assume that the residue field
$\scrO_E/\scrP_E$ has $q_E$ elements.  We denote by $|\cdot|_E$
the absolute value $E$, such that $|\varpi_E|_E=q_E\inv$, for any generator
$\varpi_E$ of $\scrP_E$.

Now let $E$ be a local non-archimedean field, of characteristic zero
with an involution $\theta$.  In certain situations
it will be more convenient to denote $\theta$
by conjugation, so, as a matter of notation, we set
\startdisp
\overline{e}=\theta(e),\;\text{for all}\; e\in E.
\finishdisp  Let $F$ be the fixed field
of $\theta$.  Since $F$ is again a local non-archimedean
field of characteristic zero, all of the above notation
again applies to $F$.  Further, we may write
\starteqn\label{eqn:omegadefn}
E=F\oplus F\omega,\;\text{as a vector space, where}\;
\omega\in E-F,\;\omega^2=a\in F-\{0\}.
\finisheqn
See \S\ref{subsec:quadraticextgeneral}, for a proof.  Note that we have
\starteqn\label{eqn:omegabar}
\theta(\omega):=\overline{\omega}=-\omega,
\finisheqn
and the relation \eqref{eqn:omegabar} completely determines
the involution $\theta$ of $E$.

Let $\scrV$ be a pair $(V,b)$ consisting
of an $m$-dimensional vector space $V$ over $E$ and a sesqui-linear
form $b$ on $V$ such that
\startdisp
\theta(b(v,u))= b(u,v)\;\text{for all}\; u,v\in V.
\finishdisp
We will also assume that $b$ is non-degenerate.

Unless  otherwise mentioned we will always denote by $G$
the group of isometries $\Isom(\scrV)$ of $\scrV$,
considered as an algebraic group over $F$.  A concrete
way of doing this is via the ``restriction of scalars"
construction.  That is we consider $\scrV$
to be a $2m$-dimensional vector space over $F$,
and then considering $G$ to be the $F$-linear
transformations of $\scrV$ satisfying an additional set of conditions
corresponding to $E$-linearity and unitarity with
respect to $\theta$.  This point of view
will be developed in greater detail when we need it,
in the proof of Lemma \ref{lem:cayleychofvar} below.

It will be convenient to fix a basis $B$ of $V$ as follows.
We fix an orthogonal $E$-basis of $\scrV$,
$B=\{v_1,\ldots, v_m\}$, such
that
\starteqn\label{eqn:diagonalbasis}
\begin{gathered}
|b(v_1,v_1)|=\cdots =|b(v_k,v_k)|=q,\\
|b(v_{k+1},v_{k+1})|=\cdots =|b(v_m,v_m)|=1.
\end{gathered}
\finisheqn
The choice of an orthogonal $E$-basis for $\scrV$
satisfying \eqref{eqn:diagonalbasis} is possible
by Th\'{e}or\`{e}me \Roman{romnine}.6.1.1 of Bourbaki
\textit{Alg\`{e}bre}, \cite{bourbakialgebra}.

For $x\in\mathrm{Mat}_n(E)$, let
\startdisp
x^*={}^t\theta(x),
\finishdisp
where the superscripted $t$ on the left
indicates the usual transpose of the matrix, and
$\theta(x)$ denotes the ``conjugation" operation
$\theta$ applied entry-wise to $x$.
We set
\starteqn\label{eqn:Tdefn}
T=\mathrm{diag}(b(v_1,v_1),\ldots,b(v_m,v_m)).
\finisheqn
Note that
\starteqn
\label{eqn:Thermitian} T=T^*\;\text{and}\; T\inv=(T\inv)^*,
\finisheqn
since $T$ is diagonal with entries in $F$.

Using the basis $B$, we write $G=U_m(E)$ as a matrix group.
We have an isomorphism of $G$ with the group
\starteqn\label{eqn:Umdefn}
U_m(T)=\{g\in \mathrm{GL}_m(E)\;|\; g^*Tg=T\}\equiv G
\finisheqn
We will write the Lie algebra $\Gg$ of $G$ in the matrix
form
\starteqn\label{eqn:unitaryliealgdefn}
\Gg\cong\Gu_m(F):=\{x\in M_m(E)\;|\; x^* T+Tx=0\}.
\finisheqn
All representations $\pi$ of $G$, considered here, are assumed
to be admissible.  We denote by $V_{\pi}$ a vector space
realization of $\pi$, and, if it has a central character,
we denote it by $\omega_{\pi}$.  Note that the center $Z(G)$
of $G$ is isomorphic to $U_1(E)$, the elements of norm $1$ in $E$,
and $Z(G)$'s isomorphic image in $U_m(T)$ is $U_1(E)I_m$.

\section{The doubling method}
Consider $\scrV\times \scrV$, consisting
of the doubled space $V\times V$ equipped with the bilinear
form $b^*=b\oplus (-b)$.  Denote by $H$ the isometry
group of $(V\times V,b^*)$.  Since the subspace
\startdisp
V^{\triangle}=\{(v,v)\;|\; v\in V\}
\finishdisp
is an $m$-dimensional isotropic subspace of $V\times V$,
hence a maximal isotropic subspace $V\times V$, the group
$H$ is quasi-split.  The elements $(g_1,g_2)$ of $G\times G$
act on $V\times V$ by
\startdisp
(g_1,g_2)(v_1,v_2)=(g_1(v_1),g_2(v_2)),
\finishdisp
and they clearly preserve $b^*$.  Thus we get a natural embedding
$i: G\times G\hookrightarrow H$.  Consider the maximal parabolic
subgroup $P_{V^\triangle}$ of $H$ which preserves $V^\triangle$.
This is a Siegel type parabolic subgroup of $H$.  Its Levi
part is isomorphic to $\mathrm{GL}(V^\triangle)\cong \mathrm{GL}(V)$.
Denote the unipotent radical of $P_{V^\triangle}$ by $U_{V^\triangle}$.
We have the ``transversality".
\starteqn\label{eqn:parabolicembeddedcopyinters}
i(G\times G)\cap P_{V^\triangle}=i(G^\triangle),
\finisheqn
where
\startdisp
G^\triangle=\{(g,g)\;|\; g\in G\}.
\finishdisp
Recall that $P_{V^\triangle}\backslash H/i(G\times G)$ is finite
and contains only one open orbit, which is $P_{V^\triangle}
\cdot i(G\times G)=P_{V^\triangle}\cdot (G\times 1)$.
This equality follows from \eqref{eqn:parabolicembeddedcopyinters}, as in
\cite{gelbartpsrallis85}, p. 8.  Denote by $\det(\cdot)$
the algebraic character of $P_{V^\triangle}$ given by
$P\mapsto\det (P|_{V^\triangle})$.
Let $\chi$ be a (unitary) character of $E^*$.
Consider, for $s\in\Complex$,
\startdisp
\rho_{\chi, s}=\mathrm{Ind}_{P_V^\triangle}^H(\chi\circ\det\cdot)
|\det\cdot|^{s-\half}.
\finishdisp
The induction is normalized as in \S 3 of \cite{rallissoudry05}.

Let $\pi$ be an irreducible representation of $G$, acting
in a space $V_{\pi}$.  Consider the contragredient
representation $\hat{\pi}$ acting in $\hat{V}_{\pi}$,
the smooth dual of $V_{\pi}$.  Denote by $\langle\cdot,\cdot\rangle$
the canonical $G$-invariant bilinear form on $V_{\pi}\times
\hat{V}_{\pi}$.  Let $v_1\in V_{\pi}$ and $\hat{v}_2\in\hat{V}_{\pi}$,
and let $f_{\chi,s}\in V_{\rho_{\chi,s}}$ be a holomorphic
section.  The local zeta integrals attached to $\pi$ by the doubling
method are
\starteqn\label{eqn:zetaintegral}
Z(v_1,\hat{v}_2,f_{\chi, s})=\int_G\langle\pi(g) v_1,
\hat{v}_2\rangle f_{\chi, s}(i(g,1))\intd g.
\finisheqn
By Theorem 3 in \cite{lapidrallis}, the integral
in \eqref{eqn:zetaintegral} converges absolutely in a right-half
plane and continues to a meromorphic function in the whole
plane.  This function is rational in $q^{-s}$.  We keep
denoting the analytic continuation by $Z(v_1,\hat{v}_2,f_{\chi, s})$.
Consider the intertwining operator
\startdisp
M(\chi, s)=M(s): \rho_{\chi, s}\rightarrow \rho_{\theta(\chi)\inv, 1-s},
\finishdisp
defined, first for $\Rept(s)\gg 0$, as an absolutely convergent integral
\starteqn\label{eqn:intertwiningfactordefn}
M(\chi, s)f_{\chi, s}(h)=\int_{U_{V^\triangle}}f_{\chi, s}(wuh)\intd u,
\finisheqn
and then, by meromorphic continuation, to the whole plane.
Here, we take, as in \cite{lapidrallis}, $w=i(1,-1)\in i(G\times G)$.
Note that
\startdisp
w(V^\triangle)=V^{-\triangle}=\{(v,-v)\;|\; v\in V\}.
\finishdisp
The subspace $V^{-\triangle}$ is a maximal isotropic subspace
of $V\times V$ which is transversal to $V^\triangle$, \textit{i.e.},
$V^\triangle\cap V^{-\triangle}=\{0\}$.  By Theorem 3 in
\cite{lapidrallis}, we have
a functional equation (as an identity of meromorphic functions
in the whole plane)
\starteqn\label{eqn:functionaleq}
\Gamma(\pi,\chi,s)Z(v_1,\hat{v}_2,f_{\chi,s})=Z(v_1,\hat{v}_2,
M(\chi,s)f_{\chi,s})
\finisheqn
for all $v_1\in V_{\pi}$, $\hat{v}_2\in \hat{V}_{\pi}$, $f_{\chi, s}\in
V_{\rho, s}$ (holomorphic section).  The function
$\Gamma(\pi, \chi, s)$ depends on the choice of measure $\intd u$
made in the definition of $M(\chi, s)$.
\vspace*{0.3cm}

Section 9 of \cite{lapidrallis} explains how to obtain
the local gamma factor $\gamma(\pi,\chi,s,\psi)$ from $\Gamma(\pi,\chi,s)$.
In the Hermitian case, the relation between the two is
given in (25), \cite{lapidrallis} as
\starteqn\label{eqn:littleandbiggammafactor}
\gamma(\pi,\chi,s,\psi)=\frac{\xi_{G}(\chi,A)}{C_H(\chi,s,A,\psi)}
\omega_{\pi}(-1)\Gamma(\pi,\chi,s),
\finisheqn
where $C_H(\chi,s,A,\psi)$ is a certain rational function of $q^{-s}$,
which depends only on $\chi$, $\psi$, a certain matrix A, and H---see
\S 5 of \cite{lapidrallis} for the exact definition---
and where $\xi_G(\chi,A)=\chi\inv(\det A)|\det A|^{s-\half}$.  It
is shown in \cite{lapidrallis}, \S\S 8--9, that $\gamma(\pi,\chi,s,\psi)$
is independent of the choice of $A$.

Since by \eqref{eqn:littleandbiggammafactor}, $\gamma(\pi,\chi,s,\psi)$
is obtained from $\omega_{\pi}(-1)\Gamma(\pi,\chi,s)$ by multiplication
by a factor which depends only on $\chi,\psi$ (and $G$), and
not on $\pi$, Theorem \ref{thm:generalstability} will follow from
the explicit formula for $\omega_{\pi}(-1)\Gamma(\pi,\chi,s)$
in Theorem \ref{thm:main},
valid for $\chi$ sufficiently ramified, which evidently does not depend
on $\pi$.

\noindent
\begin{thm} \label{thm:main}  Let $\pi$ be an irreducible
representation of $G$.  Then there exists a positive integer
$N$, such that for any ramified character $\chi$ of $E^{\times}$,
having conductor $1+\scrP_E^{N_{\chi}}$ with $N_{\chi}>N$,
we have
\startdisp
\omega_{\pi}(-1)\Gamma(\pi,\chi,s)=|2|_F^{m^2}\chi^{-m}(-1)
\times\hspace*{-.30cm}
\int\limits_{\Gg(\scrP_E^{n_{\chi}})}\hspace*{-.25cm}
\chi\inv(\det(I_m-v))
\intd \mu(v),
\finishdisp
where
\startdisp
n_{\chi}:=\left[\frac{N_{\chi}+1}{2}\right].
\finishdisp
The measure $\intd\mu(v)$ is to be specified below,
in \eqref{eqn:dmuydefn}.
\end{thm}
\vspace*{0.3cm}

We end this section with an explicit description of $i(g,1)$,
$g\in G$ as a matrix, following the decomposition
\starteqn\label{eqn:vtimesvdecomp}
V\times V=
V^{\triangle}\oplus V^{-\triangle}
\finisheqn
and a choice of a \textbf{standard basis} of $V\times V$ whose
Gram matrix, with respect to $b^*$ is $w_{2m}$.
Here, we are using the notation
\startdisp
w_n=\begin{pmatrix}&&1\\&\cdot&\\1&&\end{pmatrix},\;
\text{for}\; n\in \Natural,
\finishdisp
as in \cite{rallissoudry05}.
Let $T$ be the diagonal matrix representing $b$, as defined in
\eqref{eqn:Tdefn}.  In order to obtain a standard basis of $V\times V$,
in the sense of a basis consistent with the
decomposition \eqref{eqn:vtimesvdecomp}
whose Gram matrix (with respect to $b^*$) is $w_{2m}$, we proceed
as follows.  We let $u_i=\frac{1}{b(v_i,v_i)}v_i$, for $i=1,\ldots, m$.
Then $\{u_1,\ldots, u_m\}$ is a basis of $V$, dual to $B$
(with respect to $b$).  Then
\starteqn\label{eqn:standardbasis}
\tilde{B}=\left\{(v_1,v_1),\ldots, (v_m,v_m),\half(u_m,-u_m),
\ldots, \half(u_1,-u_1)\right\}
\finisheqn
is a standard basis of $V\times V$.
Writing the elements of $G$ as matrices, with respect to $B$,
and the elements of $H$ as matrices, with respect to $\tilde{B}$,
it is now easy to verify Lemma \ref{lem:standardbasisexpression}.
\begin{lem} \label{lem:standardbasisexpression}
We have, for all $g\in G$,
\starteqn\label{eqn:standardbasisexpression}
i(g,1)=\begin{pmatrix}\half(g+I_m)&\frac{1}{4}(g-I_m)T\inv w_m\\
                       w_mT(g-I_m)&\half w_m T(g+I_m)T\inv w_m
\end{pmatrix}
\finisheqn
and
\starteqn\label{eqn:wstandardbasisexpression}
w=i(1,-1)=\begin{pmatrix}& \half T\inv w_m\\
                    2w_mT&                \end{pmatrix}
\finisheqn
\end{lem}
\section{Proof of Theorem 3.1}\label{sec:mainproof}
Choose $f_{\chi, s}$ so that it is supported in the open orbit
$P_{V^{\triangle}}\cdot i(G\times G)=P_{V^{\triangle}}\cdot i(G\times 1)$,
and so that the restriction $\left.f_{\chi, s}\right|_{i(G\times 1)}$,
thought of as a function of $G$, is the characteristic function
$\phi_U$ of a small compact open subgroup $U$ of $G$.
We assume that $U$ is small enough, so that $V_{\pi}^U\neq 0$.
Let $0\neq v_1\in V_{\pi}^U$, and choose $\hat{v}_2\in\hat{V}_{\pi}$,
such that $\langle v_1,\hat{v}_2 \rangle=1$.  Then the integral
\eqref{eqn:zetaintegral} converges for all $s$, and is easily seen
to be
\starteqn\label{eqn:denomofgammafactor}
Z(v_1,\hat{v}_2,f_{\chi, s})=\int_U\langle \pi(u)v_1,\hat{v}_2 \rangle
\intd u=m(U)\langle v_1,\hat{v}_2\rangle=m(U),
\finisheqn
where $m(U)$ is the measure of $U$, and so, from the functional
equation \eqref{eqn:functionaleq},
\starteqn\label{eqn:gammafactorinitialexpression}
\Gamma(\pi, \chi, s)=\frac{1}{m(U)}\int_G\langle\pi(g)v_1,\hat{v}_2\rangle
M(\chi, s)f_{\chi,s}(i(g,1))\intd g,
\finisheqn
for $\Rept(s)\ll 0$.  Our next task is to compute $M(\chi,s)f_{\chi,s}
(i(g,1))$ for our choice of $f_{\chi,s}$.

From
now till \eqref{eqn:constancyonUappliedchofvars},
we assume that $\Rept(s)\gg0$, so
that the expression for $M(\chi,s)f_{\chi,s}$ given in
\eqref{eqn:intertwiningfactordefn} is valid.
\begin{lem}\label{lem:mainlem}  For the above choice of $f_{\chi,s}$,
we have, for $\Rept(s)\gg0$,
\vspace*{-0.2cm}

\begin{multline}\label{eqn:intertwiningopevalonspecialf}
M(\chi,s)f_{\chi,s}(i(g,1))=\\|2|_F^{m(1-2s)}\chi^{-m}(-2)\times \int_{G}
\chi(\det(I_m+h))|\det(I_m+h)|_E^{s+\frac{m}{2}-\half}
\phi_U(-hg)\intd h.
\end{multline}
\vspace*{0.1cm}
\end{lem}

\begin{myproof}{Proof}
Let $\Rept(s)\gg0$, so that the integral in
\eqref{eqn:intertwiningfactordefn}
converges absolutely.  Using \eqref{eqn:standardbasisexpression},
we have
\starteqn\label{eqn:intertwiningstbasis}
\begin{aligned}
M(\chi,s)f_{\chi,s}(i(g,1))&=&&\int_{(w_mx)^*=-(w_mx)}f_{\chi,s}
\left(w\begin{pmatrix}I_m&x\\&I_m\end{pmatrix}i(g,1)\right)\intd x\\
&=&&\int_{(w_mx)^*=-(w_mx)}f_{\chi,s}\left(
\begin{pmatrix}0&\half T\inv w_m\\
           2w_mT&2w_mTx\end{pmatrix}i(g,1)\right)\intd x
\end{aligned}
\finisheqn
We choose $\intd x$ to be the standard measure of matrices which
are ``skew-hermitian with respect to the second diagonal", \textit{i.e.},
\starteqn\label{eqn:dxdefn}
\intd x=\prod_{i+j<m+1}\hspace*{-3mm}\intd x_{ij}\prod_{i=1}^{m}\intd x^{(\omega)}_{ij}
\finisheqn
where $\intd x_{ij}$
is the Haar measure of $E$ which assigns the measure $1$
to $\scrO_E$, and $\intd x^{(\omega)}_{ij}$ is the $F$-invariant
measure of $F\omega$ which assigns the measure $1$ to $\scrO_{F}\omega$.

By the choice of $f_{\chi,s}$, we must have, in order for $x$
to make a nonzero contribution to the integral,
\startdisp
w\begin{pmatrix}I_m&x\\&I_m\end{pmatrix}i(g,1)\in P_{V^{\triangle}}\cdot
i(G\times 1),
\finishdisp
\textit{i.e.},
\starteqn\label{eqn:elementofsupportcondition}
w\begin{pmatrix}I_m&x\\&I_m\end{pmatrix}\in P_{V^{\triangle}}\cdot
i(G\times 1).
\finisheqn
Explicating \eqref{eqn:elementofsupportcondition}, we must solve,
for a given $x\in \mathrm{Mat}_m(E)$, such that
$(w_mx)^*=-w_mx$,
\starteqn\label{eqn:elementofsupportcondition1}
\begin{pmatrix}0&\hspace*{-2mm}\half T\inv w_m\\
           2w_mT&\hspace*{-2mm}2w_mTx\end{pmatrix}
           =
\begin{pmatrix}
E&Y\\
0&E^{\dagger}
\end{pmatrix}
\begin{pmatrix}\half(h+I_m)&\hspace*{-2mm}\frac{1}{4}(h-I_m)T\inv w_m\\
                       w_mT(h-I_m)&\hspace*{-2mm}\half w_m T(h+I_m)T\inv w_m
\end{pmatrix}.
\finisheqn
Here $E^{\dagger}=w_m(E^*)\inv w_m$ and $Y$ is such that
\startdisp
\begin{pmatrix}
E&Y\\
0&E^{\dagger}
\end{pmatrix}\in P_{V^{\triangle}}.
\finishdisp
The condition \eqref{eqn:elementofsupportcondition1} results
from \eqref{eqn:elementofsupportcondition} by substituting
the expressions of \eqref{eqn:intertwiningstbasis} and
\eqref{eqn:standardbasisexpression}.

The system \eqref{eqn:elementofsupportcondition1} has
a solution, if and only if
\starteqn\label{eqn:systemsolutioncondition}
\det(I_m-2xw_mT)\neq 0,
\finisheqn
in which case, we get
\starteqn\label{eqn:systemsolutionhE}
\begin{gathered}
h=\frac{2xw_mT+I_m}{2xw_mT-I_m}=\frac{2xw_mT+I_m}{-\conj(T\inv)
(I_m+(2xw_mT))^*}\\
E=-(2xw_mT+I_m)\inv.
\end{gathered}
\finisheqn

Note that
\starteqn\label{eqn:determinantrelations}
\overline{\det(2xw_mT+I_m)}=(-1)^m(\det(-\conj(T\inv)
(I_m+(2xw_mT))^*)).
\finisheqn

By \eqref{eqn:systemsolutionhE} and \eqref{eqn:determinantrelations},
we have
\startdisp
\det h=(-1)^m\frac{\overline{\rho}}{\rho},\;\;\text{where}
\; \rho=\det(-
(I_m+(2xw_mT))^*).
\finishdisp
Thus, $\det h(\overline{\det h})=1$, which is to say that
$h\in\det\inv(U_1(E))$, and in particular $h$ is in $\GL{m}{E}$.
Consider the Cayley transform
\startdisp
\Gc:\gl{m}{E}'\rightarrow \mathrm{GL}_m'(E),
\finishdisp
given by
\startdisp
\Gc(y)=\frac{I_m+y}{I_m-y},
\finishdisp
where
\starteqn\label{eqn:domainrestriction}
\begin{gathered}
\gl{m}{E}'=\{y\in\Gg\;|\;\det(I_m+y)(I_m-y)\neq 0\},\\
\mathrm{GL}_m(E)'=\{t\in G\;|\;\det(t+I_m)\neq 0\}.
\end{gathered}
\finisheqn
Now observe that $\Gc$ is a bijection from $\gl{m}{E}'$
to its image $\mathrm{GL}_m'(E)$.  This is most
easily proved by the following method.  We write down a
``formal" inverse for $\Gc$,
\startdisp
\Gc\inv(t)=\frac{t-I_m}{t+I_m},
\finishdisp
meaning that the left- and right-compositions
of $\Gc\inv$ with $\Gc$ are the identity mappings (formally).
Then it follows that $\Gc\inv$ is an actual
inverse to $\Gc$ on the largest ``natural domain"
for $\Gc$ (\textit{i.e.} set excluding points where $\Gc$
fails to be defined for obvious reasons) intersected
with the inverse image under $\Gc$ of the largest
natural domain of $\Gc\inv$.  It is easy
to calculate that these natural domains and their inverse
images are as described in \eqref{eqn:domainrestriction}.

The restriction of $\Gc$ from $\gl{m}{E}'$
to $\Gg':=\mathfrak{u}_m(E)\cap\gl{m}{E}'$ is a bijection onto
the image $U_m'(E):=U_m(E)\cap\GL{m}{E}'$.  Note that
in the case of the restriction to the Lie algebra of the
unitary group, we can drop explicit mention of
the requirement that $\det(I_m+y)\neq 0$,
since $\det (I_m+y)$ is merely the conjugate $\det(I_m-y)$,
and the requirement that the latter is nonzero implies
that the former is as well.

Now, \eqref{eqn:systemsolutionhE} means that
\starteqn\label{eqn:hcayleymapform}
h=-\Gc(2xw_mT).
\finisheqn
It is easy to show that
\starteqn\label{eqn:domainofcayley}
(w_mx)^*=-(w_mx)\;\text{if and only if}\; 2xw_mT\in\Gg.
\finisheqn
Here the Lie algebra $\Gg=\Gu_m(F)$ is written in matrix form as in
\eqref{eqn:unitaryliealgdefn}.
 From \eqref{eqn:intertwiningstbasis},
\eqref{eqn:elementofsupportcondition1},
\eqref{eqn:systemsolutionhE},
\eqref{eqn:determinantrelations}, and \eqref{eqn:hcayleymapform},
we get, for $\Rept(s)\gg0$,
\begin{multline}\label{eqn:intertwiningstbasis1}
M(\chi,s)f_{\chi,s}(i(g,1))\\
=|2|_F^{-m^2}\chi(-1)^m\times\\
\times\int\limits_{\stackrel{(w_mx)^*=-(w_mx)}{\det(I_m-xw_mT)\neq 0}}
\hspace*{-10mm}\chi\inv(\det(I_m-xw_mT))
|\det(I_m-xw_mT)|_E^{-s-m/2+1/2}\phi_U(-\Gc(xw_mT)g)
\intd x.
\end{multline}

Letting $y=xw_mT$, we see by \eqref{eqn:domainofcayley} that the
domain of integration in \eqref{eqn:intertwiningstbasis1} in the variable
$y$ is $\Gg'$.  Denote by
\starteqn\label{eqn:dmuydefn}\text{
$\intd\mu(y)$ the measure $\intd\mu(y)=\intd x$, where $\intd x$
is as in \eqref{eqn:dxdefn}.}
\finisheqn
Then
\begin{multline}\label{eqn:intertwiningstbasis2}
M(\chi,s)f_{\chi,s}(i(g,1))\\
=|2|_F^{-m^2}\chi(-1)^m\times\\
\times\int\limits_{\Gg'}
\chi\inv(\det(I_m-y))
|\det(I_m-y)|_E^{-s-m/2+1/2}\phi_U(-\Gc(y)g)
\intd\mu(y).
\end{multline}
Now, in \eqref{eqn:intertwiningstbasis2},
we want to make the change of variable
\starteqn\label{eqn:cayleymap}
\Gc(y)=h.
\finisheqn
\begin{lem}\label{lem:cayleychofvar}
Let $f$ be a function in $C_c^{\infty}(G)$ such that the function
\startdisp
y\mapsto |f(\Gc(y))||\det(I_m-y)|_E^{-m},
\finishdisp
defined on $\Gg'=\Gu'$,
and extended by $0$ to $\Gg$, is integrable.  Then there
is a choice of Haar measure on $G$ such that
\starteqn\label{eqn:cayleychofvar}
\begin{aligned}
\int\limits_{\Gg'}f(\Gc(y))|\det(I_m-y)|_E^{-m}\intd\mu(y)&=&&
\int\limits_{G}f(h)|\det(h)|_E^{m/2}\intd h\\
&=&&\int\limits_{G}f(h)\intd h.
\end{aligned}
\finisheqn
\end{lem}
\begin{myproof}{Proof}
Since the Jacobian of \eqref{eqn:cayleymap} is given by a rational
function, defined over $F$, it is enough to compute it over an
algebraic closure $\overline{F}$ of $F$ containing $E$.
Thus, we may assume, for this proof, that $F$ is algebraically
closed.
Over the algebraically closed field $F$, $G$ becomes the split
group $\GL{m}{F}$.  See \S\ref{subsec:formandsplitgp}
for the elementary verification of this last statement.
  Throughout this proof we will use $|\cdot|$
to denote $|\cdot|_{F}$.  We may replace
the first equality of \eqref{eqn:cayleychofvar},
to be proved, with the new formula
\starteqn\label{eqn:cayleychofvaralgclosed}
\int\limits_{\gl{m}{F}'}f(\Gc(y))|\det(I_m-y)|^{-2m}\intd\mu(y)=
\int\limits_{\GL{m}{F}'}f(h)|\det(h)|^m\intd h.
\finisheqn
We compare the two sides of
\eqref{eqn:cayleychofvaralgclosed} by computing
both using the appropriate forms of Weyl's Integration Formula.
Because, over the field $F$, which by the above argument may
be assumed to be algebraically closed, there is only one
conjugacy class of Cartan subgroups of $\GL{m}{F}$, the
``algebra" form of Weyl's Integration Formula
says that the left-hand side of \eqref{eqn:cayleychofvaralgclosed}
is
\startmultline\label{eqn:weylformulaalg}
\int\limits_{\scrL(F)}|D(x)|\int\limits_{\GL{m}{F}/L(F)}f(\Gc(\conj(g)x))|
\det(I_m-\conj(g)x)|^{-2m}\intd g\intd x=\\
=\int\limits_{\scrL(F)}|D(x)|\int\limits_{\GL{m}{F}/L(F)}f(\conj(g)(\Gc(x)))|
\det(I_m-x)|^{-2m}\intd g\intd x=\\
=\int\limits_{\scrL(F)}|D(x)||\det(I_m-x)|^{-2m}\hspace*{-.3cm}
\int\limits_{\GL{m}{F}/L(F)}f(\conj(g)(\Gc(x)))
\intd g\intd x.
\end{multline}
\vspace*{0.3cm}

\noindent
Here, $\conj$ is the conjugation map, $\scrL(F)$ is
the diagonal subalgebra of
$\gl{m}{F}$, $L$ is the diagonal subgroup of $\GL{m}{F}$; $D(x)$
is the coefficient of $t^m$ in the characteristic polynomial
$\det(tI_{\mathfrak{gl}_m}-\ad(x))$ of the linear endomorphism $\ad(x)$.
Write
\startdisp
x=\mathrm{diag}(x_1,\ldots, x_m),\;\,\text{with}\; x_i\in F.
\finishdisp
And denote by
\startdisp
\intd x=\intd x_1\cdots\intd x_m
\finishdisp
the product measure.
Since $D(x)$ is the product of the \textit{roots} of $\GL{m}{F}$,
evaluated at $x$, we have
\starteqn\label{eqn:capDcomputation}
|D(x)|=\prod_{1\leq i<j\leq m}|x_i-x_j|^2
\finisheqn
Change variables in \eqref{eqn:weylformulaalg}, $\frac{1+x_i}{1-x_i}=t_i$,
$i=1,\ldots, n$ (reflecting the explicit formula for $\Gc$).
Then $x_i=\frac{t_i-1}{t_i+1}$ (reflecting the explicit
formula for $\Gc\inv$), and $\intd x_i=|2|\frac{|t_i|}{|t_i+1|^2}\intd^*t_i$.
Using this and \eqref{eqn:capDcomputation}, a simple calculation shows that
\eqref{eqn:weylformulaalg} equals, up to a positive constant
(a power of $|2|$),
\starteqn\label{eqn:weylformulaalgt}
\int\limits_{L(F)}\prod_{1\leq i<j\leq m}
|t_i-t_j|^2\prod_{i=1}^m|t_i|
\left(\int\limits_{\GL{m}{F}/L(F)}
f(\conj(g)t)\intd g\right)
\intd^* t
\finisheqn
Here, we use the notation
\startdisp
\begin{gathered}
t=\mathrm{diag}(t_1,\ldots, t_m),\;\text{for}\;t_i\in F-\{0\},
\;i=1,\ldots, m,\\
d^*t=d^*t_1\cdots d^*t_m
\end{gathered}
\finishdisp

The ``group" form of the Weyl integration formula on $\GL{m}{F}$
says that the right-hand side of \eqref{eqn:cayleychofvaralgclosed}
is, up to a positive constant,
\startdisp
\int\limits_{L(F)}|d(t)|\left(\int\limits_{\GL{m}{F}/L(F)}
f(\conj(g)t)|\det(\conj(g)t)|^m\intd g\right)\intd^* t,
\finishdisp
where $d(t)$ is the coefficient of $z^m$ in the polynomial
$\det(zI_{\mathfrak{gl}_m}-(\Ad(t)-\mathfrak{gl}_m))$.
By clearly $\det(\conj(g)t)=\det(t)$, and the resulting
factor of $|\det(t)|^m$ in the inner integrand can
be taken out of the inner integral, and we obtain, for
the right-hand side of \eqref{eqn:cayleychofvaralgclosed},
\starteqn\label{eqn:weylformulagroup}
\int\limits_{L(F)}|d(t)||\det(t)|^m\left(\int\limits_{\GL{m}{F}/L(F)}
f(\conj(g)t)\intd g\right)\intd^* t,
\finisheqn
Since $d(t)$ is the product over all the roots of $\GL{m}{F}$
of the difference of the root from $1$, evaluated at $t$,
we have
\starteqn\label{eqn:smallDcomputation}
|d(t)|=\prod_{1\leq i<j\leq m}|t_i-t_j|^2\prod_{i=1}^m|t_i|^{-m+1}.
\finisheqn
Since
\startdisp
|\det(t)|=\prod_{i=1}^m |t_i|,
\finishdisp
the equality \eqref{eqn:smallDcomputation}
implies that \eqref{eqn:weylformulagroup} equals
\starteqn\label{eqn:weylformulagroupdeval}
\int\limits_{L(F)}\prod_{1\leq i<j\leq m}|t_i-t_j|^2
\prod_{i=1}^m |t_i|
\hspace*{-.1cm}\left(\int\limits_{\GL{m}{F}/L(F)}
f(\conj(g)t)\intd g\right)\intd^* t,
\finisheqn
Since \eqref{eqn:weylformulaalgt} and \eqref{eqn:weylformulagroupdeval}
are equal, the first equality of Lemma \ref{lem:cayleychofvar} is proved.

For the justification that $|\det(h)|=1$ for $h\in G$, see
\S\ref{subsec:fieldnormextensions}.
\end{myproof}

We continue with the proof of Lemma \ref{lem:mainlem}.  We make
the change of variable \eqref{eqn:cayleymap} in
\eqref{eqn:intertwiningstbasis2} ($\Rept(s)$ is still large enough),
and by \eqref{eqn:cayleychofvar}, we get that there is a choice
of Haar measure $\intd h$ on $G$ such that
\begin{multline*}
M(\chi,s)f_{\chi,s}(i(g,1))\\
=|2|_F^{m(1-2s)}\chi^{-m}(-2)\times \int_{G}
\chi(\det(I_m+h))|\det(I_m+h)|_E^{s+\frac{m}{2}-\half}
\phi_U(-hg)\intd h.
\end{multline*}
This proves Lemma \ref{lem:mainlem}.
\end{myproof}

The following facts about the Cayley transform are easy
to verify.
\begin{lem} \label{lem:cayleyprops} Let $v\in\Gg'$ and $g\in G'$.
Assume that the elements $I_m\pm v\in \mathrm{Mat}_m(E)$ are invertible.  Then
\begin{itemize}
\item[\textbf{CT 1}] $I_m-\Gc(v)g=(I_m-v)\inv(-\Gc\inv(g)-v)(I_m+g)$.
\item[\textbf{CT 2}] $I_m+\Gc(v)g=(I_m-v)\inv(I_m+v\Gc\inv(g))(I_m+g)$.
\item[\textbf{CT 3}] We have $I_m+\Gc(v)g$ invertible if and only
if $I_m+v\Gc\inv(g)$ is invertible, and in this case,
examining the ratio of \textbf{CT 1} to \textbf{CT 2},
\startdisp
\Gc\inv(\Gc(v)g)=(I_m-v)\inv(\Gc\inv(g)+v)(I_m+v\Gc\inv(g))\inv(I_m-v).
\finishdisp
\item[\textbf{CT $1'$}] $I_m-\Gc(v)\inv g=(I_m+v)\inv
(-\Gc\inv(g)+v)(I_m+g)$.
\item[\textbf{CT $2'$}]  $I_m+\Gc(v)\inv g=(I_m+v)\inv(I_m-v\Gc\inv(g))
(I_m+g)$.
\item[\textbf{CT $3'$}]  We have $I_m+\Gc(v)\inv g$ is invertible
if and only if $I_m-v\Gc\inv(g)$ is invertible, and in that case,
examining the ratio of \textbf{CT $1'$} to \textbf{CT $2'$},
\startdisp
\Gc\inv(\Gc(v)\inv g)=(I_m+v)\inv(-\Gc\inv(g)+v)(I_m-v\Gc\inv(g))\inv
(I_m+v).
\finishdisp
\end{itemize}
\end{lem}

Denote, for a matrix $x\in \mathrm{Mat}_m(E)$,
\starteqn\label{eqn:matrixnormdefn}
||x||=\max_{1\leq i,j\leq m}|x_{ij}|_{E}.
\finisheqn
This is a norm on $\mathrm{Mat}_m(E)$.  It is not difficult to verify
that the norm
satisfies these five properties.
\begin{itemize}
\item[\textbf{Norm 1}]  $||x+y||\leq \max\{||x||,||y||\},\;\text{for all}\;
x,y\in\mathrm{Mat}_m(E).$
\item[\textbf{Norm 2}]  $||x+y||=||y||,\; \text{if}\; ||x||<||y||.$
\item[\textbf{Norm 3}]  $||k_1xk_2||=||x||,\;\text{for all}\;
k_1,k_2\in\GL{m}{\scrO_E}.$
\item[\textbf{Norm 4}]  If $||v||<1$, then we have
\startdisp
I_m-v\in\GL{m}{\scrO_E},
\finishdisp
meaning that $I_m-v$ is both invertible and integer (possessed
of integral entries).
\item[\textbf{Norm 5}]  $||x^*||=||x||$.
\end{itemize}
\vspace*{0.2cm}
For the justification of \textbf{Norm 5}, see Lemma
\ref{lem:normandabsval}(c).
Let $N$ be a positive even integer such that
\starteqn\label{eqn:Nqcondition}
q^N>|8|_E \inv q^4,
\finisheqn
and such that
\starteqn\label{eqn:NUcontainmentcondition}
\Gc(\Gg(\scrP_E^{\frac{N-2}{2}}))\subset U.
\finisheqn
From now on, we assume that the conductor $1+\scrP_E^{N_{\chi}}$
of $\chi$ is such that $N_\chi>N$. We now return to the integral
in \eqref{eqn:intertwiningopevalonspecialf}, evaluated at
$-g\inv$ in place of $g$.  In order for a point $h$
in the domain of integration to make a non-zero contribution
to the integral (\textit{i.e.}, in order for the integrand
to be nonzero at $h$), we must have $hg\inv=u$, where $u\in U$.
That is, we must have $h=ug$ for some $u\in U$.
Then, according to \eqref{eqn:intertwiningopevalonspecialf}, we
have
\begin{multline}\label{intertwiningopevaluchofvars}
M(\chi,s)f_{\chi,s}(i(-g\inv,1))=|2|_F^{m(1-2s)}\chi^{-m}(-2)\times
\sum_{L=-\infty}^{\infty}\\
\int\limits_{\stackrel{||\Gc\inv(ug)||=q^L}{u\in U}}
\chi(\det(I_m+ug))|\det(I_m+ug)|_E^{s+\frac{m}{2}-\half}
\intd u,
\end{multline}
\vspace*{0.3cm}

\noindent
for $\Rept(s)\gg 0$.  Denote, for $g\in G$, $\Rept(s)\gg 0$,
and $L\in\Int$,
\starteqn\label{eqn:ILintegraldefn}
I_L(\chi,s; g) :=\int\limits_{\stackrel{||\Gc\inv(ug)||=q^L}{u\in U}}
\chi(\det(I_m+ug))|\det(I_m+ug)|_E^{s+\frac{m}{2}-\half}
\intd u.
\finisheqn
Set
\starteqn\label{eqn:littlenchidefn}
n_{\chi}=\left[\frac{N_{\chi}+1}{2}\right].
\finisheqn
Our main aim from this point is to prove an analogue of Lemmas 4.4--6
from \cite{rallissoudry05}, namely that
\begin{lem} \label{lem:largeL} We have
\startdisp
I_L(\chi,s; g)=0,
\finishdisp
for all $L\geq -n_{\chi}$.
\end{lem}
The first main step towards proving Lemma \ref{lem:largeL}
will be making a change of variable in the integral
$I_L(\chi,s; g)$ that will allow us to replace the single
integral of $I_L(\chi,s; g)$ with a double integral.
In order to state Lemma \ref{lem:singleintegraltodouble}, it
is convenient to introduce the following piece of notation.
\begin{defn}  Let $\preplus(\cdot)$ be the ``non-negativity
function" from the reals to the non-negative reals.
That is, let $\preplus(\cdot)$ be defined
piecewise by
\startdisp
\preplus r=\begin{cases}r&\text{if $r\geq 0$}\\
0&\text{if $r<0$}
\end{cases}.
\finishdisp
Let $\preminus(\cdot)$ be the ``non-positivity function"
from the reals to the non-negative reals defined analogously so that
\starteqn\label{eqn:nonpositivityfunctiondefn}
\text{for all $r\in\Real$},\;r=\preplus r-\preminus r,\;
\text{and}\; |r|=\preplus r+\preminus r.
\finisheqn
\end{defn}
\begin{lem}\label{lem:singleintegraltodouble}
Suppose that $L$ satisfies
\starteqn\label{eqn:Lgreaternchi}
L>-n_{\chi},\;\text{equivalently}\; \preminus L<n_{\chi}.
\finisheqn
Then
\begin{multline}\label{eqn:doubleintegralexpforIL}
I_L(\chi,s; g)=\frac{1}{\mu(\Gg_{{}^+\hspace*{-.7mm}L +n_{\chi}})}
\int\limits_{\Gg_{\preplus L+n_{\chi}}}\\
\int\limits_{\stackrel{||\Gc\inv(ug)||=q^L}{u\in U}}
\chi(\det(I_m+\Gc(v)ug))|\det(I_m+\Gc(v)ug)|_E^{s+\frac{m}{2}-\half}
\intd u
\intd\mu(v).
\end{multline}
\end{lem}

We start with a lemma of a preliminary nature, useful in making
changes of variable.
\begin{lem} \label{lem:chofvarsaid} Let $\Gg$, $||\cdot||$,
be as above, $L\in\Int$, $U$
any subgroup of $G$ and $\intd  u$ a left Haar measure on $U$.
For $A\subset G$, let $\mathbf{1}_A$ denote the characteristic
function of $A$.  Let $a,g\in G$ be fixed, and assume that these
elements of $G$
satisfy the following two conditions,
\starteqn\label{eqn:aeltofU}
a\in U,
\finisheqn
and
\starteqn\label{eqn:chofvarsaidcondition}
\text{for}\; u\in U,\;
||\Gc\inv(aug)||=L\;\text{if and only if}\;||\Gc\inv(ug)||=L.
\finisheqn
Let $f$ be a function on $G$.
Assuming that either integral converges, we have
\starteqn\label{eqn:chofvarsaid}
\int_{U}f(u)\mathbf{1}_{\left(||\Gc\inv(\cdot g)||\right)\inv(L)}(u)\intd u=
\int_{U}f(au)\mathbf{1}_{\left(||\Gc\inv(\cdot g)||\right)\inv(L)}
(u)\intd u.
\finisheqn
In other words
\starteqn\label{eqn:chofvarsaidprime}
\int\limits_{\stackrel{||\Gc\inv(ug)||=q^L}{u\in U}}f(u)\intd u=
\int\limits_{\stackrel{||\Gc\inv(ug)||=q^L}{u\in U}}f(au)\intd u.
\finisheqn
\end{lem}
\begin{myproof}{Proof}
Each side of \eqref{eqn:chofvarsaid} is equal to
\startdisp
\int_{U}f(au)\mathbf{1}_{\left(||\Gc\inv(\cdot g)||\right)\inv(L)}(au)
\intd u,
\finishdisp
the left-hand side by the condition of \eqref{eqn:aeltofU}
that $a\in U$ and because
$\intd  u$ is a Haar measure on $U$.  The right-hand side
is equal to the above expression because \eqref{eqn:chofvarsaidcondition}
means that
\startdisp
U\cap \left(||\Gc(\cdot g) ||\right)\inv(L)=U\cap
\left(||\Gc\inv(a\cdot g)||\right)\inv(L).
\finishdisp
Since we are integrating over $U$, we can drop the intersections
with $U$ for the purposes of the current argument and substitute
\startdisp
\left(||\Gc\inv(a\cdot g)||\right)\inv(L)=
\left(||\Gc\inv(\cdot g)||\right)\inv(L),
\finishdisp
on the right-hand side of \eqref{eqn:chofvarsaid}.  Then,
we use the obvious equality
\startdisp
\mathbf{1}_{\left(||\Gc\inv(a\cdot g)||\right)\inv(L)}(u)=
\mathbf{1}_{\left(||\Gc\inv(\cdot g)||\right)\inv(L)}(au)
\finishdisp
to complete the proof of \eqref{eqn:chofvarsaid}.

We derive \eqref{eqn:chofvarsaidprime}
from \eqref{eqn:chofvarsaid}  by reinterpreting the multiplication
of the integral
by the characteristic function as a restriction of the domain
of integration.
\end{myproof}

\begin{myproof}{Proof of Lemma \ref{lem:singleintegraltodouble}}
Let
\starteqn\label{eqn:vdefn}
v\in\Gg_{\preplus L+n_{\chi}}
\left(:=\Gg(\scrP_E^{\preplus L+n_{\chi}})\right),\;
\text{but otherwise
arbitrary}.
\finisheqn
We can apply Lemma \ref{lem:chofvarsaid}
with $a=\Gc(v)$ to rewrite $I_L$, provided that we verify
\eqref{eqn:aeltofU} and
\eqref{eqn:chofvarsaidcondition} in the present
instance, which is the aim of the following argument.

As for \eqref{eqn:aeltofU}, it is easy to see that
\starteqn\label{eqn:cayleyofvinU}
\Gc(v)\in U,
\finisheqn
because, $\preplus L\geq 0$ (by definition),  which together with
\eqref{eqn:littlenchidefn},
and the assumption on the conductor $N_{\chi}$ just following
\eqref{eqn:NUcontainmentcondition}, implies
\startdisp
\preplus L+n_{\chi}\geq n_{\chi}=\left[\frac{N_{\chi}+1}{2}\right]
>\left[\frac{N+1}{2}\right]>\frac{N-2}{2},
\finishdisp
and then we can apply \eqref{eqn:NUcontainmentcondition}
to conclude \eqref{eqn:cayleyofvinU}.

The verification of \eqref{eqn:chofvarsaidcondition}
consists in verifying two implications.  In both
directions, we will use the fact that \eqref{eqn:vdefn}, combined with
\textbf{Norm 4}, implies that
\starteqn\label{eqn:idmminusvinvertintegrality}
I_m-v\in\GL{m}{\scrO_E}.
\finisheqn
Note that we will also use \eqref{eqn:Lgreaternchi} in the proof
of both implications.
\vspace*{0.3cm}

\noindent
\textit{First implication of \eqref{eqn:chofvarsaidcondition}.}
\hspace*{0.5mm}  Suppose $||\Gc\inv(ug)||=q^L$.

We first claim that $I_m+\Gc(v)ug$ is invertible.
Lemma \ref{lem:cayleyprops} can be applied because
of \eqref{eqn:idmminusvinvertintegrality}.  Applying
property \textbf{CT 3} from Lemma \ref{lem:cayleyprops} with
$``g"=ug$, we have
\starteqn\label{eqn:invertibilityequiv}
I_m+\Gc(v)ug\;\text{is invertible if and only if}\; I_m+v\Gc\inv(ug)\;
\text{is invertible}.
\finisheqn
However, it can be shown directly that $I_m+v\Gc\inv(ug)$ is invertible
by using the definitions of $v$ in \eqref{eqn:vdefn} and $I_L(\chi,s;g)$
in \eqref{eqn:ILintegraldefn} to observe that
\starteqn\label{eqn:difffromidmsmall}
||v\Gc\inv(ug)||\leq q^{-(\preplus L+n_{\chi})}||\Gc\inv(ug)||=
q^{-n_{\chi}+(L-\preplus L)}=q^{-n\chi-\preminus L}<1,
\finisheqn
where we have used \eqref{eqn:nonpositivityfunctiondefn}
in the last equality.
Therefore, Property \textbf{Norm 4} implies that
\starteqn\label{eqn:Imminusshiftedug}
I_m+v\Gc\inv(ug)\in\GL{m}{\scrO_E}\;\text{and in particular, invertible}.
\finisheqn
We conclude from \eqref{eqn:Imminusshiftedug} and
\eqref{eqn:invertibilityequiv} that $I_m+\Gc(v)ug$ is invertible.

Since $I_m+\Gc(v)ug$ is invertible, property \textbf{CT 3}
from Lemma \ref{lem:cayleyprops} gives an expression
for $\Gc\inv(\Gc(v)ug)$.  Using this expression and some previously
established facts, we deduce that
\starteqn\label{eqn:shifthasnoeffect}
\begin{aligned}
||\Gc\inv(\Gc(v)ug)||&=&&||(I_m-v)\inv(\Gc\inv(ug)+v)(I_m+v\Gc\inv(ug))
(I_m-v)||\\
&=&&||\Gc\inv(ug)+v||\;\text{(by \eqref{eqn:idmminusvinvertintegrality},
\eqref{eqn:Imminusshiftedug}, and \textbf{Norm 3})}\\
&=&&||\Gc\inv(ug)||,\;\text{(by \eqref{eqn:vdefn} and
\textbf{Norm 2})}\\
&=&&q^{L}\;(\text{by assumption}).
\end{aligned}
\finisheqn

\noindent
\textit{Second implication of \eqref{eqn:chofvarsaidcondition}.}
\hspace*{0.5mm}  Suppose $||\Gc\inv(\Gc(v)ug)||=q^L$.

This is similar to the first direction, but with \textbf{CT $1'$}
through \textbf{CT $3'$} playing the role of \textbf{CT 1}
through \textbf{CT 3}.  The main steps in the proof are as follows.
First, we have
\startdisp
I_m+ug\;\text{is invertible if an only if}\;I_m-v\Gc\inv(\Gc(v)ug)\;
\text{is invertible.}
\finishdisp
By similar arguments to those used in the proof of the first
implication, we show that
\startdisp
I_m-v\Gc\inv(\Gc(v)ug)\in\GL{m}{\scrO_E}.
\finishdisp
Property \textbf{CT $3'$} now applies to give
\starteqn\label{eqn:shifthasnoeffectreverse}
\begin{aligned}
||\Gc\inv(ug)||&=&&||\Gc\inv\Gc(v)ug||\\
&=&&||(I_m+v)\inv(-Gc\inv(\Gc(v)ug)+v)(I_m-v\Gc\inv(\Gc(v)ug))\inv
(I_m+v) ||\\
&=&&||-\Gc\inv(\Gc(v)ug)+v ||\\
&=&&||\Gc\inv(\Gc(v)ug) ||\\
&=&& q^L.
\end{aligned}
\finisheqn
The penultimate step is justified by the hypothesis that
$||\Gc\inv(\Gc(v)ug)||=q^L=q^{\preplus L-\preminus L}>q^{-n\chi+\preplus L}$,
(by \eqref{eqn:Lgreaternchi} again)
which is clearly at least as large as $||v||=q^{-n_{\chi}-\preplus L}$,
since $\preplus L\geq 0$.  We then apply property \textbf{Norm 2}
to obtain the penultimate equality of \eqref{eqn:shifthasnoeffectreverse}.

Combining \eqref{eqn:cayleyofvinU}, \eqref{eqn:shifthasnoeffect},
and \eqref{eqn:shifthasnoeffectreverse},
we see that the element
$a=\Gc(v)$ of $G$ satisfies
the two hypotheses of Lemma \ref{lem:chofvarsaid} .  Applying
\eqref{eqn:chofvarsaidprime}, with
\startdisp
f(u)=\chi(\det(I_m+ug))|\det(I_m+ug)|_E^{s+\frac{m}{2}-\half},
\finishdisp
and using the definition of $I_{L}(\chi,s; g)$ in \eqref{eqn:ILintegraldefn},
we deduce that for $v$ as in \eqref{eqn:vdefn},

\begin{multline*}
I_L(\chi,s; g,v):=\int\limits_{\stackrel{||\Gc\inv(ug)||=q^L}{u\in U}}
\chi(\det(I_m+\Gc(v)ug))|\det(I_m+\Gc(v)ug)|_E^{s+\frac{m}{2}-\half}
\intd u=\\
I_L(\chi,s; g).
\end{multline*}
\vspace*{0.3cm}

\noindent
In other words, the value of $I_L(\chi,s; g,v)$, defined
in the first line of the above equality, is independent of
$v$ satisfying \eqref{eqn:vdefn}
and equal to $I_L(\chi,s; g,v)$.  Obviously, averaging
the constant $I_L(\chi,s; g,v)$
over $\Gg_{\preplus L+n_{\chi}}$ and dividing by the total measure
$\mu(\Gg_{\preplus L+n_{\chi}})$ does nothing, so we have
\startdisp
I_L(\chi,s; g)=\frac{1}{\mu(\Gg_{\preplus L+n_{\chi}})}
\int\limits_{\Gg_{\preplus L+n_{\chi}}}
I_L(\chi,s; g,v)\intd\mu(v),
\finishdisp
where $\intd\mu(v)$ is any invariant measure on $\Gg$.  In particular,
with $\intd\mu(v)$ the measure on $\Gg$ in \eqref{eqn:dmuydefn},
we complete the derivation of
\eqref{eqn:doubleintegralexpforIL} by using
the definition of $I_L(\chi,s; g,v)$ in the previous expression
for $I_L(\chi,s; g)$.
\vspace*{0.3cm}
\end{myproof}

\begin{lem} \label{lem:doubletoiteratedintegral}
Assume that \eqref{eqn:Lgreaternchi} is satisfied.
Continuing with the calculation of $I_L$
from Lemma \ref{lem:singleintegraltodouble}, there exist
$a, b\in E^{\times}$, satisfying the properties
\starteqn\label{eqn:aandbnorms}
|a|_E=|b|_E=q^{N_{\chi}},
\finisheqn
and
\starteqn\label{eqn:aplusbnorm}
|a+b|_E\leq q^{n_{\chi}},
\finisheqn
\hspace*{0.3cm}
such that
\begin{multline*}
I_L=\frac{1}{\mu(\Gg_{\preplus L+n_{\chi}})}\int
\limits_{\stackrel{||\Gc\inv(ug)||=q^{L}}{u\in U}}
\chi(\det(I_m+ug))|\det(I_m+ug)|^{s+\frac{m}{2}-\half}\times\\
\times \int\limits_{v\in\Gg_{\preplus L+n_{\chi}}}\psi_0(
b\tr v +a\tr
(v\Gc\inv(ug)))\intd \mu(v)\intd u,
\end{multline*}
\hspace*{0.3cm}

\noindent
where $\psi_0$
is a fixed character of $E$ whose conductor
is $\scrO_E$.
\end{lem}
\begin{myproof}{Proof}
By \eqref{eqn:vdefn} and \eqref{eqn:difffromidmsmall} we can
apply Lemma \ref{lem:absvalsanddets}, parts (b) and (c),
respectively, to conclude that
\starteqn\label{eqn:Imminusvdeterminantabsval}
|\det(I_m-v)|\;\text{and}\; |\det(I_m+v\Gc\inv(ug))|=1.
\finisheqn

Using Property \textbf{CT 2}, we may rewrite the integrand
of \eqref{eqn:doubleintegralexpforIL} as
\begin{multline}\label{eqn:doubleintegrandreduced}
\chi\inv(\det((I_m-v))\chi(\det(I_m+v\Gc\inv(ug)))\chi(\det(I_m+ug))\times\\
\times |\det(I_m-v)|_E^{-s+\frac{m}{2}+\half}
|\det(I_m+v\Gc\inv(ug))|_E^{s-\frac{m}{2}-\half}
|\det(I_m+ug)|_E^{s-\frac{m}{2}-\half}=\\
=\chi\inv(\det((I_m-v))\chi(\det(I_m+ug))\chi(\det(I_m+v\Gc\inv(ug)))
|\det(I_m+ug)|_E^{s-\frac{m}{2}-\half}.
\end{multline}
\vspace*{0.3cm}

\noindent
By using
\eqref{eqn:Imminusvdeterminantabsval}, we see that the
fourth and fifth factors in the product don't contribute
to the first expression, to obtain the latter expression
of \eqref{eqn:doubleintegrandreduced}.

We now apply Corollary \ref{cor:characterrewriting}
to rewrite each of the factors
$\chi\inv(\det((I_m-v))$
and $\chi(\det(I_m+v\Gc\inv(ug)))$.  In the first case
$n=\preplus L+n_{\chi}$, while in the second
$n=n_{\chi}$, and in both cases $N=2n$.  Since $N_{\chi}\leq 2n$,
the hypothesis \eqref{eqn:condchibounds} is verified.  The Corollary
now gives $a,b\in E^{\times}$ satisfying \eqref{eqn:aandbnorms},
with
\starteqn\label{eqn:characterrewriting}
\begin{aligned}
\chi\inv(\det((I_m-v))\chi(\det(I_m+v\Gc\inv(ug)))&=&&
\psi_0(b\tr v)\psi_0(a\tr(v\Gc\inv(ug)))\\
&=&&\psi_0(b\tr v+a\tr(v\Gc\inv(ug)))
\end{aligned}
\finisheqn
for all $v$ such that $v,v\Gc\inv(ug)\in\scrP_E^{n_{\chi}}$.
By \eqref{eqn:difffromidmsmall} \eqref{eqn:characterrewriting} does
 indeed apply to the $v$ considered in the integrand of
\eqref{eqn:doubleintegralexpforIL}.
Together with the calculation of \eqref{eqn:doubleintegrandreduced},
this proves the Lemma with the exception of the claim
\eqref{eqn:aplusbnorm}.

In order to verify \eqref{eqn:aplusbnorm}, note that,
by construction, $a,b$ are elements of $E^\times$
such that
\startdisp
\chi(1+x)=\psi_0(ax),\;\text{and}\;\chi\inv(1+x)=\psi_0(bx),\,\;
\text{for all}\; x\in\scrP_E^{n_{\chi}}.
\finishdisp
Consequently, we have
\startdisp
q^0=1=\psi_0(ax+bx)=\psi((a+b)x)\;\text{for all}\; x\in\scrP_E^{n_{\chi}},
\finishdisp
implying that
\startdisp
(a+b)\scrP^{n_{\chi}}_E\subseteq \scrO_E,\;\text{so that}\;
|a+b|_Eq^{-n_{\chi}}\leq 1.
\finishdisp
This immediately yields \eqref{eqn:aplusbnorm}.
\end{myproof}

Rewriting the single integral
$I_L(\chi,s;g)$ as the iterated integral of \eqref{eqn:doubleintegralexpforIL}
was the first main step in showing that $I_L(\chi,s;g)=0$.
From this point, the strategy of the proof consists in changing
the order of integration and rewriting the inner integrand in such
a way that this inner integrand can be seen to be zero,
for all $u$ in the range of integration of the outer integrand.
\begin{lem}\label{lem:innerintegralvanishing}
Let $L\in\Int$ and \bmth\textbf{assume that $L\neq 0$}\ubmth.
Let $\psi_0$
be a fixed character of $E$ with conductor $\scrO_E$,
$a,b\in E^{\times}$ such that
\starteqn\label{eqn:aabsval}
|a|_E=|b|_E>|2|_E\inv q^{n\chi+1}.
\finisheqn
Then for $X\in\mathfrak{gl}_m(E)$ such that
\starteqn\label{eqn:innterintlemcond}
||X||=q^L,
\finisheqn
we have
\starteqn\label{eqn:innerintegralvanishing}
\int\limits_{\Gg_{\preplus L + n_{\chi}}}\psi_0(b\tr v+a\tr(vX))
\intd\mu(v)=0.
\finisheqn
\end{lem}
\begin{myproof}{Proof}
Note that the expression $b\tr v+a\tr(vX)$ is equal to
\startdisp
\tr(v(b I_m+aX))
\finishdisp
As in the proof of Lemma 4.4 of \cite{rallissoudry05}, at (4.42)
we see that in order for the integral
\eqref{eqn:innerintegralvanishing} to be nonzero, we must have
\starteqn\label{eqn:conditionfornonzerointeg}
||bI_m+aX||\leq |2|_E\inv q^{\preplus L+n_{\chi}+1}.
\finisheqn
Since $L\neq 0$, by assumption, $q^L\neq 1$, so that
\eqref{eqn:innterintlemcond} and the equality of
\eqref{eqn:aabsval} together imply that
\startdisp
||bI_m||\neq ||aX||.
\finishdisp
Therefore, Property \textbf{Norm 2} implies that
\starteqn\label{eqn:norm2appliedtosum}
||bI_m+aX||=\max(||bI_m||,||aX||).
\finisheqn
By the strict inequality of \eqref{eqn:aabsval},
\eqref{eqn:innterintlemcond}, and the definition
of $\preplus L$ before \eqref{eqn:nonpositivityfunctiondefn}
\startdisp
\max(||bI_m||,||aX||)> |2|_E\inv q^{\preplus L+n_{\chi}+1},
\finishdisp
so that by \eqref{eqn:norm2appliedtosum}
\startdisp
||bI_m+aX||> |2|_E\inv q^{\preplus L+n_{\chi}+1}.
\finishdisp
With \eqref{eqn:conditionfornonzerointeg} this gives a contradiction.
Therefore, by the above comments, the integral
of \eqref{eqn:innerintegralvanishing} equals zero.
\end{myproof}

\begin{lem}
\label{lem:innerintegralvanishingLzero}
Let $\psi_0$
be a fixed character of $E$ with conductor $\scrO_E$,
$a,b\in E^{\times}$ satisfying \eqref{eqn:aandbnorms} and
\eqref{eqn:aplusbnorm}.
Then for
\starteqn\label{eqn:XinGzero}
X\in\Gg,\;\text{such that}\;\, ||X||=1,
\finisheqn
we have
\starteqn\label{eqn:innerintegralvanishingLzero}
\int\limits_{\Gg_{n_{\chi}}}\psi_0(b\tr v+a\tr(vX))
\intd\mu(v)=0.
\finisheqn
\end{lem}
\begin{myproof}{Proof}
In order for the integral of \eqref{eqn:innerintegralvanishingLzero} not
to vanish, we must have
\starteqn\label{eqn:conditionfornonzerointegLzero}
||bI_m+aX||\leq |2|_E\inv q^{n_{\chi}+1},
\finisheqn
paralleling \eqref{eqn:conditionfornonzerointeg}.  The hypothesis
\eqref{eqn:aplusbnorm} is equivalent to assuming
that $a+b$ is of the form $\varpi^{-n_{\chi}}\Go$, for $\varpi$
the generator of $\scrP_E$ and some $\Go\in \scrO_E$.  So
in order for the integral of \eqref{eqn:innerintegralvanishingLzero}
to be nonzero, we must have
\startdisp
||\varpi^{-n_{\chi}}\Go+a(X-I_m)||\leq |2|_E\inv q^{n_{\chi}+1},
\finishdisp
By \textbf{Norm 2}, because $||\varpi^{-n_{\chi}}\Go||\leq q^{n_{\chi}}$,
this implies that
\startdisp
||a(X-I_m)||\leq |2|_E\inv q^{n_{\chi}+1},
\finishdisp
so by \eqref{eqn:aandbnorms},
\startdisp
||X-I_m||\leq |2|_E\inv q^{n_{\chi}-N_{\chi}+1},
\finishdisp
meaning that
\startdisp
X-I_m\in\mathfrak{gl}_m(\scrP_E^{N_{\chi}-n_{\chi}-1-\nu_E(2)}).
\finishdisp
Thus $X\in I_m+\mathfrak{gl}_m(P^{N_{\chi}-n-1-\nu_E(2)})$.  Because
of \eqref{eqn:Nqcondition}, this, combined with
Proposition \ref{prop:stayingawayfrom1}, implies that $X\notin\Gg$.
Thus, we obtain a contradiction with \eqref{eqn:XinGzero}.
\end{myproof}

\begin{myproof}{Completion of Proof of Lemma \ref{lem:largeL}}
By applying Lemma \ref{lem:innerintegralvanishing}
or \ref{lem:innerintegralvanishingLzero} to
Lemma \ref{lem:doubletoiteratedintegral}, depending on whether
$L=0$ or $L\neq 0$, we deduce
the vanishing of $I_L$ for all $L$ in the required
range of $L\geq -n_{\chi}$.  The reason the hypotheses
of Lemma \ref{lem:innerintegralvanishing} are satisfied
is that according to Lemma \ref{lem:doubletoiteratedintegral}
we have \eqref{eqn:aandbnorms} and then \eqref{eqn:Nqcondition}
implies \eqref{eqn:aabsval}.  The reason the hypotheses
of Lemma \ref{lem:innerintegralvanishingLzero} are satisfied in the case
$L=0$ is that $X=\Gc\inv(ug)$ satisfies
\eqref{eqn:XinGzero} for all $u$ in the domain
of integration of the outer integral for $I_L$.
\end{myproof}

\begin{myproof}{Completion of Proof of Theorem \ref{thm:main}}
Combining \eqref{eqn:ILintegraldefn} and \eqref{eqn:littlenchidefn}
and Lemma \ref{lem:largeL},
we have, for $\Rept(s)\gg 0$,
\begin{multline}\label{eqn:ILvanishingapplied}
M(\chi,s)f_{\chi,s}(i(-g\inv,1))\\=|2|_F^{m(1-2s)}\chi^{-m}(-2)
\times\hspace*{-0.88cm}
\int\limits_{\stackrel{||\Gc\inv(ug)||\leq
q^{-n_{\chi}}}{u\in U}}\hspace*{-0.75cm}
\chi(\det(I_m+ug))|\det(I_m+ug)|_E^{s+\frac{m}{2}-\half}
\intd u.
\end{multline}
\vspace*{0.3cm}

\noindent
The domain of integration in \eqref{eqn:ILvanishingapplied} is over
$u\in U$, such that
\startdisp
ug\in\Gc(\Gg(\scrP_E^{n_{\chi}}))\subset\Gc\left(\Gg
\left(\scrP_E^{\frac{N-2}{2}}\right)\right)
\subseteq U,
\finishdisp
where the latter two containments follow from \eqref{eqn:NUcontainmentcondition}
and \eqref{eqn:littlenchidefn}.  This implies that $M(\chi,s)f_{\chi,s}
(i(-g\inv,1))$ is supported in $U$.  Since we may now assume that $g\in U$
and $\intd u$ is the Haar measure on $U$ with total mass $1$,
the following general
observation applies.  For any function $f$ on $U$
we have
\startdisp
F\;
\text{defined by, for $g\in U$},\;
F(g):=\int_{U} f(ug)\intd u,\;\text{is a constant function
equal to $F(1)$}.
\finishdisp
We apply this general observation with $f$ the product of
the integrand in \eqref{eqn:ILvanishingapplied} and the characteristic
function of the the set
\startdisp
\Gc\inv\left(\left(-\infty,q^{-n_{\chi}}\right)\right).
\finishdisp
The result of doing so is to deduce that $M(\chi,s)f_{\chi,s}
(i(-g\inv,1))$ is constant on $U$ and equal to
$M(\chi,s)f_{\chi,s}(i(-1,1))$.  Note further that
\startdisp
|\det(I_m+u)|=|2^m\det(I_m-(I_m-u)/2)|=|2|^m\left|\det\left(I_m-
\frac{I_m-u}{2}\right)\right|
\finishdisp
and the latter factor is $1$ for $u$ close to $I_m$,
by Lemma \ref{lem:absvalsanddets}, part (b).  Therefore,
the entire factor $|\det(I_m+u)|_E^{s+\frac{m}{2}-\half}$
in the integrand of \eqref{eqn:ILvanishingapplied}
comes out of the integral as $|2|_E^{m(s+\frac{m}{2}-\half)}=
|2|_F^{m(2s+m-1)}$.  These arguments imply that,
that the expression for $M(\chi,s)f_{\chi,s}(i(-g\inv,1))$
in \eqref{eqn:ILvanishingapplied} can be replaced by
\starteqn\label{eqn:constancyonUapplied}
M(\chi,s)f_{\chi,s}(i(-1,1))=|2|_F^{m^2}\chi^{-m}(-2)
\times\hspace*{-.9cm}
\int\limits_{||\Gc\inv(u)||\leq
q^{-n_{\chi}}}\hspace*{-.85cm}
\chi(\det(I_m+u))
\intd u.
\finisheqn
We may change variable in \eqref{eqn:constancyonUapplied},
$u=\Gc(v)$, $v\in\Gg(\scrP_E^{n_{\chi}})$,
and choose the measure $\mu(v)$ on $\Gg(\scrP_E^{n_{\chi}})$
as in \eqref{eqn:dmuydefn}, and get (using $I_m+\Gc(v)=2(I_m-v)\inv$),
\starteqn\label{eqn:constancyonUappliedchofvars}
M(\chi,s)f_{\chi,s}(i(-1,1))\\=|2|_F^{m^2}\chi^{-m}(-1)
\times\hspace*{-.35cm}
\int\limits_{\Gg(\scrP_E^{n_{\chi}})}\hspace*{-.35cm}
\chi\inv(\det(I_m-v))
\intd \mu(v).
\finisheqn
We have proved \eqref{eqn:constancyonUappliedchofvars} for $\Rept(s)\gg 0$.
But $s$ does not appear on the right
side of \eqref{eqn:constancyonUappliedchofvars}.
Therefore, \eqref{eqn:constancyonUappliedchofvars}
is valid for all $s$ by analytic continuation.  In summary,
we have shown that, for all $s\in\Complex$,
$M(\chi,s)f_{\chi,s}(i(-g\inv,1))$ is given by the integral
of the right side of \eqref{eqn:constancyonUappliedchofvars} for $g\in U$,
independent of $g$, and by $0$ for $g\notin U$.

Now let $\Rept(s)$ be small enough that the expression of
\eqref{eqn:gammafactorinitialexpression} is valid.  By the
last sentence of the previous paragraph, we can replace the integral
over $G$ by an integral over $U$, then the factor
$M(\chi,s)f_{\chi,s}(i(-u\inv,1))$ in the integrand by
$M(\chi,s)f_{\chi,s}(i(-1,1))$.  Further, since $U$
is chosen so that $v_1$ is $U$-invariant, and further, since
$\langle v_1,\hat{v}_2\rangle=1$, the first factor
$\langle\pi(-1)v_1,\hat{v}_2\rangle$ reduces
to $\omega_{\pi}(-1)$.  So on the right side we are left integrating
the constant $M(\chi,s)f_{\chi,s}(i(-1,1))$
over $U$ the measure on $U$ was chosen to have
total mass one.  It follows then that \eqref{eqn:gammafactorinitialexpression}
reduces to
\starteqn\label{eqn:gammafactorfinalexpression}
\omega_{\pi}(-1)\Gamma(\pi,\chi,s)=M(\chi,s)f_{\chi,s}(i(-1,1)).
\finisheqn
The right-hand side of \eqref{eqn:gammafactorfinalexpression} is
independent of $\chi, s$, and depends only on $\pi$.  Therefore,
by analytic continuation, \eqref{eqn:gammafactorfinalexpression}
is valid for all $s$.  Substituting the expression of
\eqref{eqn:constancyonUappliedchofvars} for $M(\chi,s)f_{\chi,s}(i(-1,1))$
into the right-hand side of \eqref{eqn:gammafactorfinalexpression},
we obtain Theorem \ref{thm:main}.
\end{myproof}

\section{Appendix: elementary verifications}\label{sec:appendix}
\subsection{Generalities concerning quadratic extensions of local fields.}
\label{subsec:quadraticextgeneral}
We begin by proving the claim in the second paragraph of
\S\ref{sec:notandprelim}.
Let $\omega'\in E-F$.
Then, since $E$ is a quadratic extension
of $F$, there exist $a_0,b\in F$, such that
\startdisp
\omega'^2+a_0\omega'+b=0.
\finishdisp
Completing the square, we obtain the equivalent equations
\startdisp
\begin{aligned}
\omega^2+a_0\omega'+\frac{a_0^2}{4}&=&&\frac{a_0^2}{4}-b,\\
\left(\omega'+\frac{a_0}{2}\right)^2&=&&\frac{a_0^2}{4}-b.\\
\end{aligned}
\finishdisp
Now set $\omega=\omega'+a_0^2/2$ and $a=\frac{a_0}{4}-b$, to
obtain an $\omega$ and $a$ as claimed.  Note that we
have only used the field property and $\mathrm{char}E\neq 2$
in this argument (in practice $\mathrm{char}E$ is zero).
In order to verify \eqref{eqn:omegabar},
note that, since $a\in F=\mathrm{Fix}(\theta)$,
\startdisp
\overline{\omega}^2=a=\omega^2,
\finishdisp
so that we have
\startdisp
\left(\frac{\omega}{\overline{\omega}}\right)^2=1,
\finishdisp
from which we deduce
\startdisp
\overline{\omega}=-\omega,
\finishdisp
since $\omega\in E-F$, by assumption.

\bmth
\subsection{Passage from $G=U_m(E)$ to
$\GL{m}{\overline{E}}=\GL{m}{\overline{F}}$.}
\ubmth
\label{subsec:formandsplitgp}
We explicate and prove the sentence\vspace*{3mm}\newline
\hspace*{2cm}\begin{minipage}[c]{9cm}\textit{\hspace*{-1mm}
Over $\overline{F}$, defined by be an algebraic closure
of $F$ containing $E$, $G$ becomes the split group $\GL{m}{\overline{F}}$.}
\vspace*{3mm}
\end{minipage}\newline
\noindent from the proof of Lemma \ref{lem:cayleychofvar}.
In order to make the statement more comprehensible, we break
down the passage from $F$ to $\overline{F}$ into two steps, as follows.
\begin{itemize}
\item{\textbf{Step 1.}}  Replace $F$ with the quadratic extension $E$,
over which $G$ becomes $\GL{m}{E}$.
\item{\textbf{Step 2.}}  Replace $E$ with $\overline{F}=\overline{E}$,
over which $\GL{m}{E}$ becomes $\GL{m}{\overline{E}}=\GL{m}{\overline{F}}$,
hence split.
\end{itemize}
It is well-known that a split group over an algebraically
closed field has a conjugacy class of Cartan subgroups,
so this argument also validates the application of the ``simple"
form of the Weyl integration formulas, \eqref{eqn:weylformulaalg}
and \eqref{eqn:weylformulagroup}, above.

Of the two steps above, \textbf{Step 2} amounts to a straightforward
tensoring operation on the Lie algebra level.  Alternatively,
considering the abstract Chevalley group $\mathbf{GL_m}$ as a functor
from fields (\textbf{Field}) to matrix groups,
we may consider this step as a simple
substitution of objects belonging to the domain category \textbf{Field}.
Therefore, only \textbf{Step 1} needs any further explanation,
for which we introduce the following concepts.
\begin{defn}  \label{defn:restrictionofscalars}
Let $F\subset E$
be a quadratic extensions of fields of characteristic different from $2$.
By \S\ref{subsec:quadraticextgeneral}, this implies that
\textit{as a vector space},
$E$ is of the form $F\oplus\omega F$, with $\omega^2=a\in F$.
For $V$ an $m$-dimensional vector space over
$E$, let $V(F)$ be the $2m$-dimensional $F$-vector space which equals
$V$ as a set.  If
\startdisp
\beta=\{v_1,\ldots, v_m\}
\finishdisp
is a basis for $V$, then call
\startdisp
\beta^{F}:=\{v_1,\;\omega v_1,\ldots, v_m,\;\omega v_m\},
\finishdisp
the \bmth\textbf{corresponding basis for $V(F)$}.  \ubmth
Fix a $B$ as above so as to consider
$\mathfrak{gl}(V(F))$ as the matrix-algebra $\gl{2m}{F}$.
For an element $X\in\mathfrak{gl}(V(F))$, denote by $X^{(ij)}$
the $ij^{\rm th}$ \textit{two-by-two block} counting from the upper-left of the
matrix $X$, for $1\leq i,j\leq m$.
Thus $X^{(ij)}\in\gl{2}{F}$, and its entries
will be denoted by $X^{(ij)}_{kl}$ for $1\leq k,l\leq 2$.
Define the operation \texttt{transp} on $\mathfrak{gl}(V(F))$
by the relation
\startdisp
(\texttt{transp}(X))^{(ij)}=X^{(ji)},\;\text{for $1\leq i,j\leq m$},
\finishdisp
in other words, the transpose operation operating a matrix of $2$-by-$2$
blocks instead of individual entries.  Define the operation
\texttt{conj} on $\mathfrak{gl}(V(F))$ by the relations
\startdisp
(\texttt{conj}(X))^{(ij)}_{kl}=(-1)^{k-l}(X)^{(ij)}_{kl},\;
\text{for $1\leq i,j\leq m$, $1\leq k,l\leq 2$}
\finishdisp
\textit{i.e.} \texttt{conj} negates ``off-diagonal" entries
of each two-by-two block and leaves the diagonal entries alone.
Finally, let $\Gg$ be the subalgebra of $\mathfrak{gl}(V(F))$
satisfying the conditions\vspace*{0.2cm}
\begin{itemize}
\item[\bmth\textbf{$E$-linearity.}\ubmth] Each $X^{(ij)}$ is of the form
\startdisp
\begin{pmatrix}
f_1&af_2\\
f_2&f_1
\end{pmatrix}.
\finishdisp
\item[\bmth\textbf{$E$-unitarity.}
\ubmth]  We have $\texttt{conj}\circ\texttt{transp}(X)=-X$.
\end{itemize}\vspace*{0.1cm}
\end{defn}
\begin{prop}  With $\Gg$ defined as above, we have $\Gg\cong\Gu_{m}(E)$
as a Lie algebra over $F$.
\end{prop}
\begin{myproof}  {Proof} Consider the $F$-isomorphism of $V$
with $V(F)$ induced by the identification of the underlying sets.
This isomorphism induces natural embedding of $\gl{m}{E}$ into
$\gl{2m}{F}$.  The first condition, \bmth\textbf{$E$-linearity},\ubmth
\hspace*{0.5mm}
is equivalent to $X$ belonging to the image of this embedding,
that is to $X$'s actually from a ``restriction of scalars" from
an element of $\gl{m}{E}$.
Given that $X$ satisfies the condition of \bmth\textbf{$E$-linearity}\ubmth,
it is clear that the inverse image of $X$
belongs to $\Gu_m(F)$ if and only if $X$ satisfies the second condition,
of \bmth\textbf{$E$-unitarity}\ubmth.
\end{myproof}
\vspace*{0.3cm}

Therefore, in order to justify \textbf{Step 1}, we can replace
the original task of showing that $\Gu_m(F)\otimes E\cong \gl{m}{E}$,
with showing that $\Gg\otimes E\cong \gl{m}{E}$.
The isomorphism is simply the isomorphism induced on the level
of vector space endomorphisms by the vector-space isomorphism
$V(F)\otimes E\cong V$.  This can
be made obvious by showing how a basis of the former
maps to a basis of the latter.  In order to define
the basis, we adopt the notation
\startdisp
E^{(ij)}\left(\left(\begin{smallmatrix}
a&b\\c&d\end{smallmatrix}\right)\right)
\in \mathrm{Mat}_{2m}(F)\;\text{with
$\left(\begin{smallmatrix}a&b\\c&d\end{smallmatrix}\right)$ as the
$ij$-th $2$-by-$2$ block and zeros elsewhere.}
\finishdisp
Also, we use $e_{ij}$ to denote the (usual) elementary matrix
with $1$ in the $ij^{\rm th}$ position and zeros elsewhere.

Then it is easily calculated that the $E$-basis of $\Gg\otimes E$
\startdisp
\left\{E^{(ii)}\left(\begin{smallmatrix}0&a\\
1&0\end{smallmatrix}
\right)\right\}_{1\leq i\leq m}
\bigcup\left\{\left(E^{(ij)}+E^{(ji)}\right)\left(\begin{smallmatrix}&a\\
1&\end{smallmatrix}
\right)\right\}_{1\leq i<j\leq m}\bigcup\left\{
\left(E^{(ij)}-E^{ji}\right)\left(I_2
\right)
\right\}_{1\leq i<j\leq m}
\finishdisp
maps to the basis of $\gl{m}{E}$
\startdisp
\{e_{ii}\}_{1\leq i\leq m}\bigcup\,
\{\omega(e_{ij}+e_{ji})\}_{1\leq i<j\leq m}
\bigcup\, \{e_{ij}-e_{ji}\}_{1\leq i<j\leq m}.
\finishdisp

\noindent  \textbf{Example.}  As an example of the matrix form $\Gg$
for $\Gu_m$ just given, we offer the simplest case in which all
main features of the situation are visible, namely the case $m=2$.
Then we readily calculate that
\startdisp
\Gg=\left\{\left.\begin{pmatrix}&af^{(11)}_2&f^{(12)}_1&af^{(12)}_2\\
f^{(11)}_2&&f^{(12)}_2&f^{(12)}_1\\
-f^{(12)}_1&af^{(12)}_2&&af^{(22)}_2\\
f^{(12)}_2&-f^{(12)}_1&f^{(22)}_2&
\end{pmatrix}\;\right|\; f^{(ij)}_k\in F\;\text{for $i,j,k\in\{1,2\}$}\right\}.
\finishdisp
The reader can now verify by inspection the above
claims concerning
the bases
$\left\{E^{(ii)}\left(\begin{smallmatrix}0&a\\
1&0\end{smallmatrix}
\right)\right\}_{1\leq i\leq m}
\bigcup\cdots$
and
$\{e_{ii}\}_{1\leq i\leq m}\bigcup\cdots$
in this case of $m=2$ and see how to extend the arguments
to general $m$.

\bmth
\subsection{General observations concerning $\mathrm{Mat}_n(F)$,
for $F$ a local field,
and the absolute value.}
\ubmth
For this subsection, let $F$ be a non-Archimedean local field.
\begin{lem} \label{lem:absvalsanddets} Let $v\in\mathrm{Mat}_n(F)$.  Let $C\in\Int$
and assume that
\startdisp
v\in\mathrm{Mat}_n(F)_{C}:=\mathrm{Mat}_n(\scrP_F^{C}).
\finishdisp
Then we have the following,
\begin{itemize}
\item [(a)]  $\det(I_m-v)\in 1+\scrP^{C}$ and
$\det(I_m-v)\equiv 1+\tr(v)\mod \scrP^{2C}$.
\item [(b)]  If in addition $C>0$, then $|\det(I_m-v)|=1$, and
$I_m-v$ is invertible.
\item [(c)]  The conclusions of (b) hold if we make the (equivalent)
assumption that $||v||<1$, where $||\cdot||$ is the matrix norm
defined in \eqref{eqn:matrixnormdefn}.
\end{itemize}
\end{lem}
\begin{myproof}{Proof}  For (a), use induction on $m$.  The assertion
concerning the determinant is
clear for $m=1$.  For $m>1$, expand the determinant in minors
along any row or column.  All the minors
will lie in the set $1+\scrP^{C}$.  All the co-factors will lie in the ideal
$\scrP^{C}$ except for the one coming from the diagonal, which
is clearly in $1+\scrP^{C}$.  Now it easy to arrive at the conclusion.
The assertion concerning the trace is proved by an analogous induction.

For (b), use ``ultra-metric" property of the absolute value.

Part (c) follows immediately from part (b) and the definitions
of the absolute value and the matrix norm $||\cdot||$.
\end{myproof}

\begin{lem} \label{lem:groupmorphism} For $N, n$ positive integers such that
\starteqn\label{eqn:bigandlittleNrelation}
n\leq N\leq 2n,
\finisheqn
the map
\startdisp
x\mapsto 1+x
\finishdisp
defines a group isomorphism
\startdisp
\scrP_F^{n}/\scrP_F^{N}\overset{\cong}{\longrightarrow}
1+\scrP_F^n/1+\scrP_F^N.
\finishdisp
\end{lem}
\begin{myproof}{Proof}  Obviously, the map defined
in the lemma has a well-defined inverse.  Both the original map
and its inverse are defined on the whole quotient
group $\scrP_F^{n}/\scrP_F^{N}$, and
$1+\scrP_F^n/1+\scrP_F^N$, respectively.  Therefore, the
map in the lemma is a set bijection.  The condition $N\leq 2n$
is precisely what's needed to guarantee that each map respects
the group law.
\end{myproof}

It is well-known that the map
\textbf{Char} taking a group to its character group
is a (contravariant) functor on the category \textbf{Group}.
Therefore, by Lemma \ref{lem:groupmorphism},
the map
\startdisp
x\mapsto 1+x
\finishdisp
and its inverse induce an isomorphism of the character
group of $\scrP_F^{n}/\scrP_F^{N}$ with that of
$1+\scrP_F^n/1+\scrP_F^N$.  Finally, there is the well-known
description of the character group of $F$ as the mappings
of the form
\startdisp
x\mapsto \psi_0(ax),
\finishdisp
where $\psi_0$ is a fixed character of $F$ whose conductor is $\scrO_F$,
and $a\in F^*$ such that the negative of the valuation
$-v_E(a)$ equals
the conductor of $\chi$.  Summarizing
the above discussion, we have proved the following.
\begin{cor}\label{cor:characterdescription}  Let $N, n$ satisfy
\eqref{eqn:bigandlittleNrelation}.  Let
$\psi_0$ be a fixed character of $F$ whose
conductor is $\scrO_F$.  The characters
of $1+\scrP^{n}$ with conductor not exceeding $N$ are given precisely by
\startdisp
\chi(1+x)=\psi_0(ax),\;\text{for all}\; x\in\scrP_F^n,
\finishdisp
where $a$ is an element of $F^*$ satisfying $-v_F(a)$ equal
to the conductor of $\chi$.
\end{cor}

Further, putting together Corollary \ref{cor:characterdescription} and
Lemma \ref{lem:absvalsanddets}, we deduce the following
\begin{cor} \label{cor:characterrewriting} Let $n, N$ be as in
\eqref{eqn:bigandlittleNrelation},
and $\psi_0$ a fixed character of $F$ whose conductor is $\scrO_F$.
Let $\chi$ be a character of $F^{\times}$ such that
\starteqn\label{eqn:condchibounds}
\text{The conductor of $\chi$ is no greater than $N$.}
\finisheqn
Then there exists a fixed $a\in F^*$
with
\starteqn\label{eqn:valofa}
-\nu(a)\;\text{equal to the conductor of $\chi$}.
\finisheqn
such that
\starteqn
\chi(\det(1-v))=\psi_0(a\tr(v)),\;\text{for all}\;
v\in\mathrm{Mat}_m(F)\;\text{with}\;
||v||\leq q^{-n}.
\finisheqn
\end{cor}

\subsection{The Field Norm and Finite Extensions}
\label{subsec:fieldnormextensions}
Directly from the definition of $G$ in \eqref{eqn:Umdefn}, we have
\startdisp
\det(g)\overline{\det(g)}=1,
\finishdisp
so that
\startdisp
|\det(g)|_E|\overline{\det(g)}|_E=1.
\finishdisp
Since the field extension $E/F$ is quadratic, the two elements
of $\mathrm{Gal}(E/F)$ are $\Id$ and $\theta=\overline{\cdot}$.
Thus, the above equality is equivalent to
\startdisp
\bfN_{E/F}(\det (g))=1.
\finishdisp
The claim at the end of the proof of Lemma \ref{lem:cayleychofvar},
namely, that $|\det(g)|=1$, then follows from
part (b) of the following basic lemma:
\begin{lem} \label{lem:normandabsval} Let
$\Rational_p\subset F\subset E$ be a tower of finite
field extensions, with $[E:F]=m$ and $[F:\Rational_p]=n$.
Then we have
\begin{itemize}
\item[(a)]  For any $x\in E$,
\startdisp
|x|_E=\sqrt[\uproot{4}\leftroot{-4}mn]{|\bfN_{F/\Rational_p}\bfN_{E/F}(x)|_p}.
\finishdisp
\item[(b)]  The kernel of $\bfN_{E/F}$ is a subset of the
kernel of $|\cdot|_E$, in other words, for any $x\in E$,
\startdisp
\bfN_{E/F}(x)=1\;\text{implies}\; |x|_E=1.
\finishdisp
\item[(c)] $|x|_E=|\overline{x}|_E$.
\end{itemize}
\end{lem}
\begin{proof}  The proof of (a) consists of assembling several
standard facts in the theory of finite field extensions of complete
normed fields.  We use \cite{gouvea}, as a reference.  By
Corollary 5.3.2, there is at most one absolute value on $E$
extending the $p$-adic absolute value on $\Rational_p$.  By
Theorem 5.3.5, this absolute value is given by the formula
\startdisp
|x|_E=\sqrt[\uproot{4}\leftroot{-4}mn]{|\bfN_{E/\Rational_p}(x)|_{p}}.
\finishdisp
Now we rewrite the norm inside the radical using the well-known
formula
\startdisp
\bfN_{E/\Rational_p}=\bfN_{F/\Rational_p}\circ \bfN_{E/F},
\finishdisp
which appears as Problem 192 on p. 132.  This completes the proof
of (a).  Part (b) is a simple consequence of (a).

Part (c) follows directly from (a) and the calculation
\startdisp
\bfN_{E/F}(x)=x\overline{x}=\overline{x}x=\bfN_{E/F}(\overline{x}).
\finishdisp
\end{proof}

Of course, one can make much more general statements along the lines of Lemma
\ref{lem:normandabsval}, for example, replacing $\Rational_p$
with a more general non-Archimedean local field, but such
generalizations do not concern us here.

\boldmath
\subsection{Quantitative separation of $\Gg$ from $I_m$} \unboldmath In the proof of
Lemma \ref{lem:innerintegralvanishingLzero}, we needed to use
the fact that the elements of $\Gg$ cannot get ``too close"
to the identity $I_m$, in a precise quantitative sense.  The
purpose of the following proposition is simply to prove this claim.
\begin{prop}\label{prop:stayingawayfrom1}  For $T$ a form
matrix as in \eqref{eqn:diagonalbasis}--\eqref{eqn:Tdefn},
let $\Gg$ be the realization of $\mathfrak{u}_m$ given
by
\startdisp
\Gg=\{xw_mT\;|\; x\in\mathfrak{gl}_m(E)\;\text{and}\; (xw_m)^*=
-xw_m\}.
\finishdisp
Let
\startdisp
n>v_E(2).
\finishdisp
Then there is no $X\in\Gg$ such that $X-I_m\in\mathfrak{g}_m(\scrP_E^n)$.
\end{prop}
\begin{myproof}{Proof}  Suppose otherwise, so that we have
\starteqn\label{eqn:closeelement}
xw_mT-I_m\in\mathfrak{gl}_m(\scrP_E^{n}).
\finisheqn
By \eqref{eqn:diagonalbasis} and \eqref{eqn:Tdefn},
because $T\inv$ is a diagonal matrix with entries of norm $1$
or $q\inv$, hence in $\Go$.  Therefore, multiplying an element
of $\mathfrak{gl}_m(\scrP_E^{n})$ on the right by $T\inv$ multiplies
each column by an integer.
Multiply \eqref{eqn:closeelement} on the right by $T\inv$ to
obtain
\starteqn\label{eqn:gseparationlemmainter1}
xw_m-T\inv\in\mathfrak{gl}_m(\scrP_E^n).
\finisheqn
Then apply $(\cdot)^*$ to \eqref{eqn:closeelement}.
Using part (c) of Lemma \ref{lem:normandabsval}, we have
\startdisp
(xw_m)^*-(T\inv)^*\in\mathfrak{gl}_m(\scrP_E^n),
\finishdisp
so that by the description of $\Gg$ in the hypotheses and
\eqref{eqn:Thermitian}
\starteqn\label{eqn:gseparationlemmainter2}
-xw_m-(T\inv)\in\mathfrak{gl}_m(\scrP_E^n).
\finisheqn
Add \eqref{eqn:gseparationlemmainter1} and \eqref{eqn:gseparationlemmainter2}
to get
\startdisp
2T\inv\in \mathfrak{gl}_m(\scrP_E^n).
\finishdisp
By \eqref{eqn:diagonalbasis} and \eqref{eqn:Tdefn} this means
that
\startdisp
2b(v_i,v_i)\in\scrP_E^n,
\finishdisp
implying
\startdisp
|2|_Eq\;\text{or}\; |2|_E\times 1\leq q^{-n},
\finishdisp
\textit{i.e.},
\startdisp
|2|_E\leq q^{-n-1}\;\text{or}\; |2|_E\leq q^{-n}.
\finishdisp
Since $n>\nu_E(2)$, either of these condition will produce a
contradiction.
\end{myproof}
\bibliographystyle{amsalpha}
\bibliography{memoirsbib}
\end{document}

%% file: artpreamble.tex
\usepackage{graphicx}
\usepackage{amssymb}
\usepackage{relsize}
\usepackage{amsmath,amscd,verbatim}
\usepackage[mathscr]{eucal}
\newtheorem{thm}{\hspace{1cm}Theorem}[section]

\newtheorem{cor}[thm]{\hspace{1cm}Corollary}
\newtheorem{lem}[thm]{\hspace{1cm}Lemma}

\newtheorem{prop}[thm]{\hspace{1cm}Proposition}

\theoremstyle{definition}
\newtheorem{defn}[thm]{\hspace{1cm}Definition}
\theoremstyle{remark}

\numberwithin{equation}{section}

\newenvironment{myproof}[1]{\emph{#1}.  }{\qed\vspace{.3cm}}
 
 \newcommand{\Real}{\mathbb{R}}
 \newcommand{\Natural}{\mathbb{N}}
 
 \newcommand{\Complex}{\mathbb{C}}
 \newcommand{\Rational}{\mathbb{Q}}
 
 \newcommand{\Int}{\mathbb{Z}}

 \newcommand{\Gg}{\mathfrak{g}}
 \newcommand{\Gu}{\mathfrak{u}}

 \newcommand{\GL}[2]{\mathrm{GL}_{#1}(#2)}
 \newcommand{\gl}[2]{\mathfrak{gl}_{#1}(#2)}

 \newcommand{\intd}{\,\mathrm{d}}

 \newcommand{\half}{\frac{1}{2}}

 \newcommand{\Gc}{\mathfrak{c}}

 \newcommand{\ad}{\mathrm{ad}}
 
 \newcommand{\conj}{\mathbf{c}}

 \newcommand{\Ad}{\mathrm{Ad}}

 \newcommand{\Id}{\mathrm{Id}}

 \newcommand{\Go}{\mathfrak{o}}

 \newcommand{\tr}{\mathrm{tr}}
 
 \newcommand{\inv}{^{-1}}

 \newcommand{\Spos}{\mathrm{\Spos}}

 \newcommand{\Rept}{\mathrm{Re}}

 \newcommand{\finishdisp}{\vspace{0.3cm} \]}
 \newcommand{\startdisp}{\vspace{0.3cm}\[ }
 \newcommand{\starteqn}{\vspace{0.3cm}\begin{equation}}
 \newcommand{\finisheqn}{\vspace{0.3cm}\end{equation}}
  \newcommand{\startmultline}{\vspace{0.3cm}\begin{multline}}
 \newcommand{\finishmultline}{\vspace{0.3cm}\end{multline}}
 \newcommand{\starteqnarraynn}{\vspace*{0.1cm}\begin{eqnarray*}}
 \newcommand{\finisheqnarraynn}{\vspace*{2cm}\end{eqnarray*}}
  \newcommand{\starteqnarray}{\vspace*{0.1cm}\begin{eqnarray}}
 \newcommand{\finisheqnarray}{\vspace*{2cm}\end{eqnarray}}

 \newcommand{\Isom}{\textrm{Isom}}

 \newcommand{\bmth}{\boldmath}
 \newcommand{\ubmth}{\unboldmath}

 \newcommand{\scrO}{\mathscr{O}}
 \newcommand{\scrP}{\mathscr{P}}
 \newcommand{\scrV}{\mathscr{V}}
 \newcommand{\scrL}{\mathscr{L}}
 
 \newcommand{\preplus}{{}^+\hspace*{-.7mm}}
  \newcommand{\preminus}{{}^-\hspace*{-.7mm}}
  \newcommand{\bfN}{\mathbf{N}}